\documentclass[onefignum,onetabnum]{siamart220329}

\usepackage{xcolor}

\usepackage{csquotes}\usepackage[activate={true,nocompatibility},tracking=true]{microtype}\usepackage{hyphenat}\usepackage{siunitx}
\usepackage{placeins}
\usepackage{cite}\usepackage[begintext=``, endtext='']{quoting}
\SetBlockEnvironment{quoting}
\SetBlockThreshold{1}

\usepackage{amsmath}
\usepackage[fleqn]{nccmath}
\usepackage{amssymb}
\usepackage{amscd}
\usepackage{mathtools}
\usepackage[mathscr]{eucal}
\usepackage{calligra}
\DeclareMathAlphabet{\mathcalligra}{T1}{calligra}{m}{n}
\DeclareFontShape{T1}{calligra}{m}{n}{<->s*[1.1]callig15}{}
\usepackage{bm}
\usepackage{autonum}\usepackage[nice]{nicefrac}
\usepackage{multirow}

\usepackage{upgreek}
\usepackage{textgreek}

\usepackage{tikz}
\usepackage{pgf}
\usepackage{pgfplots}
\usepackage{pgfplotstable}
\usepackage{shellesc}
\pgfplotsset{compat=newest}
\usetikzlibrary{spy}
\usetikzlibrary{shapes}
\usetikzlibrary{backgrounds}
\usetikzlibrary{shadows}
\usetikzlibrary{matrix}
\usepgfplotslibrary{fillbetween}

\usepackage{comment}
\usepackage[colorinlistoftodos,prependcaption,textsize=small]{todonotes}
\usepackage{cleveref}

\crefname{equation}{}{}
\usepackage{csvsimple}
\usepackage{booktabs}
\usepackage{environ}
\usepackage{adjustbox}
\usepackage{relsize}
\usepackage{afterpage}

\usepackage[noend]{algpseudocode}

\usepackage{subcaption}

\usepackage{lipsum}
\usepackage{amsfonts}
\usepackage{graphicx}
\usepackage{epstopdf}
\usepackage{subcaption}
\usepackage{amsopn}
\ifpdf
  \DeclareGraphicsExtensions{.eps,.pdf,.png,.jpg}
\else
  \DeclareGraphicsExtensions{.eps}
\fi

\newsiamremark{remark}{Remark}
\newsiamremark{hypothesis}{Hypothesis}
\crefname{hypothesis}{Hypothesis}{Hypotheses}
\newsiamthm{claim}{Claim}

\headers{A matrix-free ILU realization}{Daniel Drzisga, Andreas Wagner, Barbara Wohlmuth}

\title{A matrix-free ILU realization based on surrogates}

\author{Daniel~Drzisga\footnotemark[2]
\and Andreas~Wagner\footnotemark[2]~\thanks{Corresponding author.
\funding{This work was partly supported by the German Research Foundation through grant WO671/11-1.}}
\and Barbara~Wohlmuth\thanks{Lehrstuhl f\"ur Numerische Mathematik, Fakult\"at f\"ur Mathematik (M2), Technische Universit\"at M\"unchen, Garching bei M\"unchen (\email{drzisga@ma.tum.de}, \email{wagneran@ma.tum.de}, \email{wohlmuth@ma.tum.de})}}

\let\min\relax \DeclareMathOperator*\min{\vphantom{p}min}
 
\let\subset\relax \DeclareMathOperator{\subset}{\subseteq}
 
\let\tilde\widetilde

\let\div\relax \DeclareMathOperator{\div}{div}

\newcommand{\sspace}{\hspace{0.25pt}}

\newcommand{\R}{\mathbb{R}}       \newcommand{\mesh}{\mcT} \newcommand{\dd}{\,\mathrm{d}}           \newcommand{\conj}[1]{\overline{#1}}

\newcommand{\btheta}{{\bm{\theta}}}

\newcommand{\mcD}{\mathcal{D}}
\newcommand{\mcE}{\mathcal{E}}
\newcommand{\mcF}{\mathcal{F}}

\newcommand{\mcI}{\mathcal{I}}

\newcommand{\mcP}{\mathcal{P}}

\newcommand{\mcS}{\mathcal{S}}
\newcommand{\mcT}{\mathcal{T}}

\newcommand{\mcV}{\mathcal{V}}

\newcommand{\bbN}{\mathbb{N}}

\newcommand{\bbP}{\mathbb{P}}

\newcommand{\bfb}{\mathbf{b}}

\newcommand{\bff}{\mathbf{f}}

\newcommand{\bfu}{\mathbf{u}}

\newcolumntype{?}{!{\vrule width 1.2pt}}

\makeatletter
\newsavebox{\measure@tikzpicture}
\NewEnviron{scaletikzpicturetowidth}[1]{	\tikzifexternalizingnext{		\def\tikz@width{#1}		\begin{lrbox}{\measure@tikzpicture}			\tikzset{external/export next=false,external/optimize=false}			\BODY
		\end{lrbox}		\pgfmathparse{#1/\wd\measure@tikzpicture}				\BODY
	}{		\BODY
	}}
\makeatother

\definecolor{color1}{rgb}{0, 0.4470, 0.7410}
\definecolor{color2}{rgb}{0.8500, 0.3250, 0.0980}
\definecolor{color3}{rgb}{0.9290, 0.6940, 0.1250}
\definecolor{color4}{rgb}{0.7060, 0.3840, 0.7650}
\definecolor{color5}{rgb}{0.4660, 0.6740, 0.1880}
\definecolor{color6}{rgb}{0.3010, 0.7450, 0.9330}
\definecolor{color7}{rgb}{0.6350, 0.0780, 0.1840}
\definecolor{color8}{rgb}{0.7410, 0.3800, 0.1840}

\newcommand{\markblack}[1]{#1}
\newcommand{\markred}[1]{\bm{#1}}
\newcommand{\markblue}[1]{{#1}}
\newcommand{\markteal}[1]{{#1}}
\newcommand{\markviolet}[1]{{#1}}

\headers{A matrix-free ILU realization based on surrogates}{Daniel Drzisga, Andreas Wagner, Barbara Wohlmuth}

\ifpdf
\hypersetup{
  pdftitle={A matrix-free ILU realization based on surrogates},
  pdfauthor={Daniel Drzisga, Andreas Wagner, Barbara Wohlmuth}
}
\fi

\begin{document}

\maketitle

\begin{abstract}
  Matrix-free techniques play an increasingly important role in large-scale simulations.
Schur complement techniques and massively parallel multigrid solvers for second-order elliptic partial differential equations can significantly benefit from reduced memory traffic and consumption.
The matrix-free approach often restricts solver components to purely local operations, for instance, to the most basic schemes like Jacobi- or Gauss--Seidel-Smoothers in multigrid methods.
An incomplete LU (ILU) decomposition cannot be calculated from local information and is therefore not amenable to an on-the-fly computation which is typically needed for matrix-free calculations.
It generally requires the storage and factorization of a sparse matrix which contradicts the low memory requirements in large scale scenarios.
In this work, we propose a matrix-free ILU realization. More precisely, we introduce a memory-efficient, matrix-free ILU(0)-Smoother component for low-order conforming finite-elements on tetrahedral hybrid grids. 
Hybrid-grids consist of an unstructured macro-mesh which is subdivided into a structured micro-mesh. 
The ILU(0) is used for degrees-of-freedom assigned to the interior of macro-tetrahedra.
This ILU(0)-Smoother can be used for the efficient matrix-free evaluation of the Steklov--Poincar\'e operator from domain-decomposition methods.
After introducing and formally defining our smoother, we investigate its performance on refined macro-tetrahedra.
Secondly, the ILU(0)-Smoother on the macro-tetrahedrons is implemented via surrogate matrix polynomials in conjunction with a fast on-the-fly evaluation scheme resulting in an efficient matrix-free algorithm.
The polynomial coefficients are obtained by solving a least-squares problem on a small part of the factorized ILU(0) matrices to stay memory efficient.
The convergence rates of this smoother with respect to the polynomial order are thoroughly studied.
\end{abstract}



\section{Introduction}\label{sec:introduction}

The incomplete LU(0)-factorization\cite{meijerink1977iterative} (ILU) approximates an LU factorization by retaining the sparsity pattern of the original matrix.
For strongly anisotropic problems in 2D, it is often used as a smoother within multigrid algorithms since its convergence rates are more stable than the ones of simpler smoothers like the Gauss--Seidel- or Jacobi-Smoothers\cite[Sec.~7.8]{wesseling1995introduction}.
This property carries on to anisotropic 3D problems in which the coupling in one spatial direction is dominant while other schemes have to be used if two of the spatial directions are dominant\cite{kettler1986aspects}. 
The related thresholded ILU-Smoother was recently used for p-multigrid in isogeometric analysis\cite{tielen2020pmultigrid,tielen2021direct} or as a smoother for the wave equation\cite{umetani2009multigrid}.

Besides its usage as a smoother, incomplete factorizations like the ILU are used as preconditioners\cite{axelsson1985incomplete,saad2003iterative}, for instance in problems involving the incompressible Stokes equation\cite{kang2020performance} or in electromagnetic scattering\cite{malas2007incomplete}. 
An algorithm for a communication avoiding ILU(0) preconditioner in the high-performance context was introduced in \cite{grigori2015communication}.
Algorithms for the efficient parallel assembly of thresholded ILU preconditioners can be found in \cite{anzt2018parilut} including adaptions to GPUs in \cite{anzt2019parilut,antz2014hybridilu}.

Matrix-free methods are becoming increasingly prevalent within finite-element frameworks\cite{kronbichler2017performance,vargas2021matrix,kohl2019hyteg}.
For instance, large scale mantle-convection simulations typically operate on scales on which storing the discretization matrices is not always feasible \cite{bauer2020terraneo}. On the other hand, reducing the memory traffic by not requiring to load a matrix from memory has the potential to result in faster algorithms on today's hardware. 
This generates interest in adapting old matrix-based algorithms to the matrix-free context.
For non-local factorization algorithms like the ILU, this poses a tremendous challenge as the matrix entries cannot be locally computed on-the-fly.

For structured grids, several techniques exist to approximate matrices for an efficient evaluation. 
For instance, stencil-scaling techniques that work for both scalar\cite{bauer2018stencil} and vectorial\cite{drzisga2020stencil} equations.
However, since the ILU-approach relies on a matrix factorization which cannot be computed locally, these approaches are inapplicable.
In our work, we propose a matrix-free ILU realization on structured subgrids based on surrogates.
Here, the discrete matrix, which usually approximates a continuous operator is additionally approximated by surrogate polynomials \cite{bauer2017two,bauer2018new,drzisga2019surrogate,bauer2019large,drzisga2020igasurrogate}.
These techniques can also be adapted for hybrid structured grids which are extensively used in \cite{bergen2005hierarchical,bergen2004hierarchical,bergen2007hierarchical,kohl2019hyteg} and consist of a coarse unstructured macro-grid which is subdivided into a fine structured micro-grid.
The former gives the approach enough flexibility to represent relevant domains while the latter provides the computational advantages of structured grids.

In this work, we apply the surrogate methodology to our factorized ILU matrix.
In the interior of the highly structured grids, we utilize an ILU factorization and approximate the resulting matrix by surrogate polynomials.
This approximation is formed in a memory-efficient way such that the memory costs stay within sensible bounds.
We therefore obtain an efficient solver in the interior of our structured grid.

To illustrate the potential of our approach, we provide two examples of how the matrix-free ILU can be used on hybrid grids:
Our main application is the approximation of the Steklov--Poincar\'e operator for the Laplacian in a matrix-free way.
This operator is a main ingredient of many non-overlapping domain-decomposition methods  and therefore efficient algorithms for its evaluation are highly relevant, see \cite{chan1994domain,korneev2004,toselli2004domain,quarteroni1999domain,mathew2008domain} and references therein. 
It formally requires the exact inversion of an elliptic equation inside a subdomain for which a multigrid method can be efficiently applied.
By using the ILU-factorization as a smoother within this inner multigrid, the inversion becomes robust with respect to distortions along one axis.
In the supplementary material a second application is provided in which we extend the subgrid ILU-Smoother to a smoother on the global grid.

The article is structured as follows:
In Section~\ref{sec:2_hybrid_grid_preliminaries}, we describe the problem, introduce the notation and present the Steklov--Poincar\'e operator.
In Section~\ref{sec:3_hybrid_ilu_smoother}, we introduce an ILU formulation that is amenable to a matrix-free algorithm.
Next, we introduce a reordering strategy on our hybrid mesh, to optimize its performance as a smoother inside single subdomains.
Finally, we introduce the matrix-free surrogate ILU in Section~\ref{sec:4_matrix_free_ilu} and compare its asymptotic convergence rates within a multigrid algorithm to the matrix-based ILU.
We conclude with a short outlook and summary in Section~\ref{sec:5_conclusion}.

\section{Hybrid grids}\label{sec:2_hybrid_grid_preliminaries}
In this section, we describe our model problem in the context of a low-order conforming finite element discretization on hybrid grids.
Hybrid grids combine the flexibility of unstructured grids with the computational advantages of structured grids \cite{bergen2005hierarchical,bergen2004hierarchical,bergen2007hierarchical,mayr2021non}. 
In addition, they provide a natural domain partitioning that can be used to distribute the work to different nodes. 
\subsection{Preliminaries and notation}

In the weak form of a Poisson-type equation, $-\div(K\nabla u) = f$, on an open domain $\Omega \subset \R^3$ with 
homogeneous Dirichlet boundary conditions on $\Gamma_D \subset \partial\Omega$, natural boundary conditions on $\partial \Omega \setminus \Gamma_D$ 
and
an inhomogeneous, bounded, symmetric, uniformly positive-definite diffusion tensor $K(x): \Omega \to \R^{3\times 3}$, we obtain the bilinear form
\begin{equation}
	a(u,v)
	=
	\int_\Omega \nabla u(x)^\top K(x)\sspace \nabla v(x) \dd x,
	\quad u, v\in V = \left\{ u \in H^1(\Omega) \,:\, u|_{\Gamma_D} = 0 \right\}
	.
\label{eq:BilinearFormPoisson}
\end{equation}
This includes the special case of a bounded, uniformly-positive scalar material parameter $\kappa: \Omega \to \R$ by setting $K = \kappa\,\textmd{Id}_3$, where $\textmd{Id}_3 \in \R^{3\times3}$ is the identity matrix.
One application of the full diffusion tensor, would be the pull-back of a blending function which maps a simple tetrahedral domain to a more complex domain, thereby providing a better approximation of the domain boundary.
Given a load $f\in L^2(\Omega)$ which defines the linear form $F(v) = \int_\Omega f\hspace{0.25pt} v \dd x$, we obtain the standard variational problem:
\(
    \text{Find } u\in V 
	\text{ satisfying }
		a(u,v) = F(v)
		\text{ for }
	v\in V 
	.
	\label{eq:PoissonStrongForm}
\)

\begin{figure}[H]
\centering
\includegraphics[width=\textwidth]{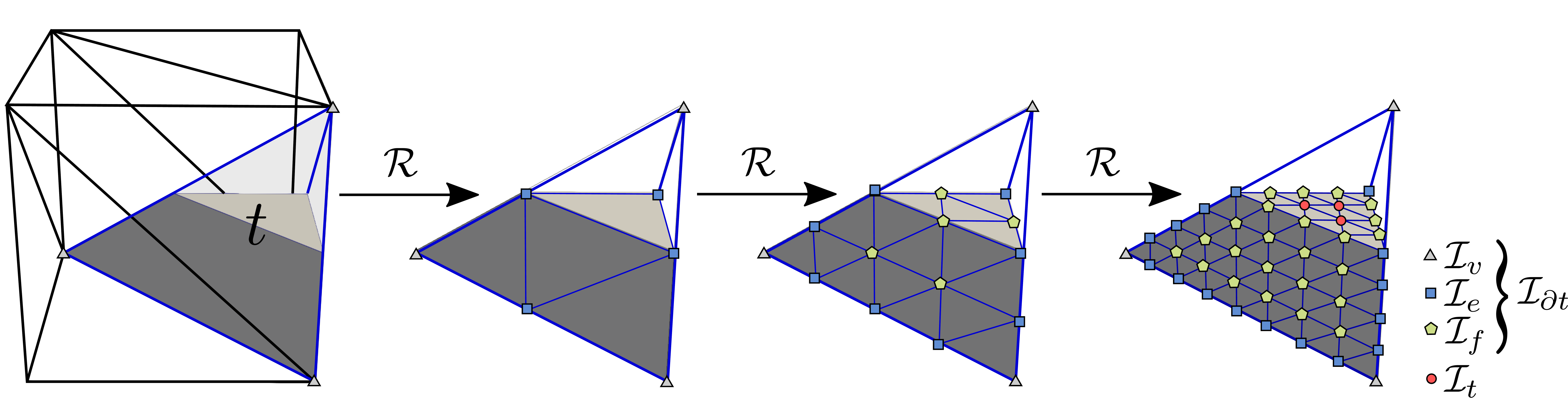}
\caption{
	Hybrid-grid refinement procedure in 3D for a clipped tetrahedron in a cubic macro-mesh.
	The DoF belonging to the index sets $\mcI_v$, $\mcI_e$, $\mcI_f$, $\mcI_t$ and $\mcI_{\partial t}$ are illustrated by different colors and shapes.
}
\label{fig:hybrid-grid}
\end{figure}

The typical approach in HHG\cite{bergen2005hierarchical,bergen2004hierarchical,bergen2007hierarchical} and HyTeG\cite{kohl2019hyteg} is 
to discretize the domain~$\Omega$ with a coarse, possibly unstructured, simplicial triangulation. 
This so-called macro-mesh consists of macro-vertices $\mathcal V_H$, macro-edges $\mathcal E_H$, macro-faces $\mathcal F_H$ and macro-tetrahedra $\mcT_H$.
All macro-primitives are referred to as $\mathcal P_H = \mathcal V_H \cup \mathcal E_H \cup \mathcal F_H \cup \mcT_H$.
Based on this initial grid, we construct a hierarchy of $L \in \mathbb{N}$, grids $\mesh = \{\mesh_{h_l},\; h_l = 2^{-l} H,\; l=2,\dots L+1\}$ by successive global uniform refinement.
The choice to start in the multigrid hierarchy with $l=2$ guarantees that each macro-element contains at least one interior element, which simplifies the notation in our algorithms.
As it is standard, each of these refinements is achieved by subdividing all elements in 3D into $8$ sub-elements.
For details of the refinement in 3D, we refer to \cite{bey1995tetrahedral}.
Due to this refinement process, the element neighborhood at each vertex in the interior of a macro element is always the same.
The whole process of the hybrid grid mesh setup is schematically depicted in Figure~\ref{fig:hybrid-grid}.

Associated with $ \mesh_{h_l}$, is the space $V_{h_l} \subset V$ of piecewise linear conforming finite elements
\[
	V_{h_l} = \{ v \in V : v|_t \in \mcP_1(t) \text{ for each } t\in \mathcal T_{h_l}\}
	.
\]

Let $\phi_i \in V_{h_l}$ and $\phi_j \in V_{h_l}$ be the scalar-valued linear nodal basis functions associated with the $i$-th and $j$-th mesh node.
The set containing all our degrees-of-freedom (DoF) indices is referred to by $\mathcal I_{h_l}$. 
If the multigrid level is obvious from the context, we will try to suppress the level dependence $h_l$ for a more compact notation. 
By $u = \sum_{i} u_i \phi_{i}$ and $v = \sum_i v_i \phi_i$ we denote linear combinations of the nodal basis function with coefficients $u_i, v_i \in \R$.
Defining the matrix $A_{ij} = a(\phi_i, \phi_j)$ and vector $\bff_i = F(\phi_i)$ results in the linear algebraic formulation of the discrete variational problem associated with the weak formulation:
\(
	\text{Find } \bfu \in \R^{|\mathcal I|}
	\text{ satisfying }
		A \bfu = \bff
	.
\label{eq:PoissonWeakFormLA}
\)

Hybrid meshes impose a domain-partitioning, which is also used for assigning the DoF in an HPC environment to computing nodes.
This approach avoids communication between the DoF located inside the same macro-primitive, while for DoF on different macro-primitives communication is necessary.
This has to be considered for an efficient evaluation of our operators since operations acting locally on the same primitive type do not require inter-node communication.
To define these local operations, we have to introduce notation to localize our vectors and matrices:
For arbitrary index sets $I \subset \mcI$ we define restriction operators $R_I: \R^{|\mcI|} \to \R^{|I|}$ consisting of zeros and ones, which discard vector entries whose component is not present in the index set and just retain entries in $I$.
We also assume that the restriction operator retains the global DoF ordering.
Given a macro-primitive $p \in \mcP_H$, we denote the set of all DoF which are located on the primitive by $\mcI_p \subset \mcI$ and its restriction operator by $R_p = R_{\mcI_p}$.
For an arbitrary macro-tetrahedron $t\in \mcT_H$, which is adjacent to macro-vertices $v_i \in \mcV_H$, $1\leq i \leq 4$,
macro-edges $e_j \in \mcE_H$, $1\leq j \leq 6$ and macro-faces $f_k \in \mcF_H$, $1\leq k \leq 4$ we define the index-set of its ghost-layer as 
$\mcI_{\partial t} = (\cup_{i=1}^4 \mcI_{v_i}) \cup (\cup_{j=1}^6 \mcI_{e_j}) \cup (\cup_{k=1}^4 \mcI_{f_k}) $.
All these sets are illustrated in Figure~\ref{fig:hybrid-grid}.

\begin{figure}[ht]
\centering
\includegraphics[width=0.25\textwidth]{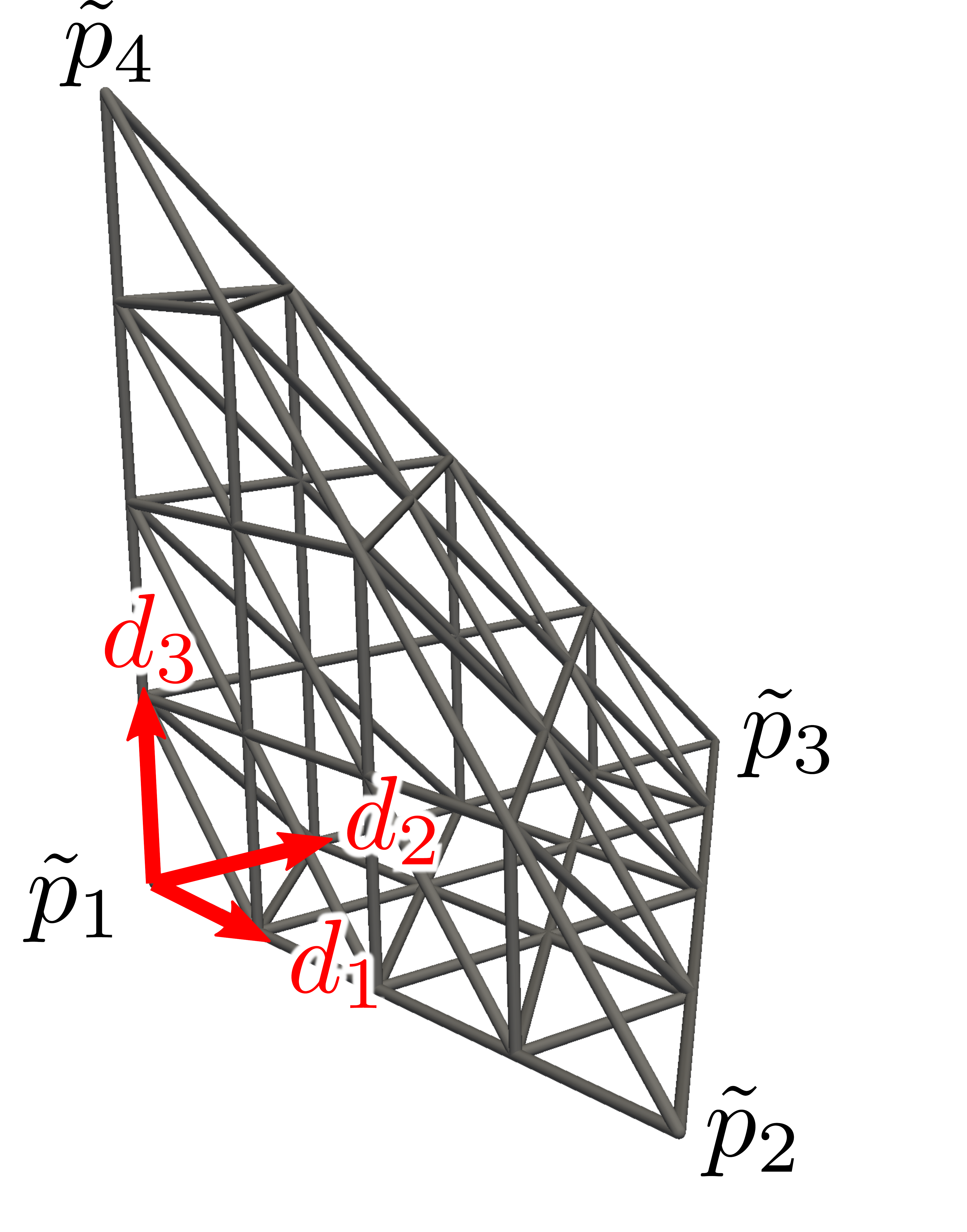}\includegraphics[width=0.50\textwidth]{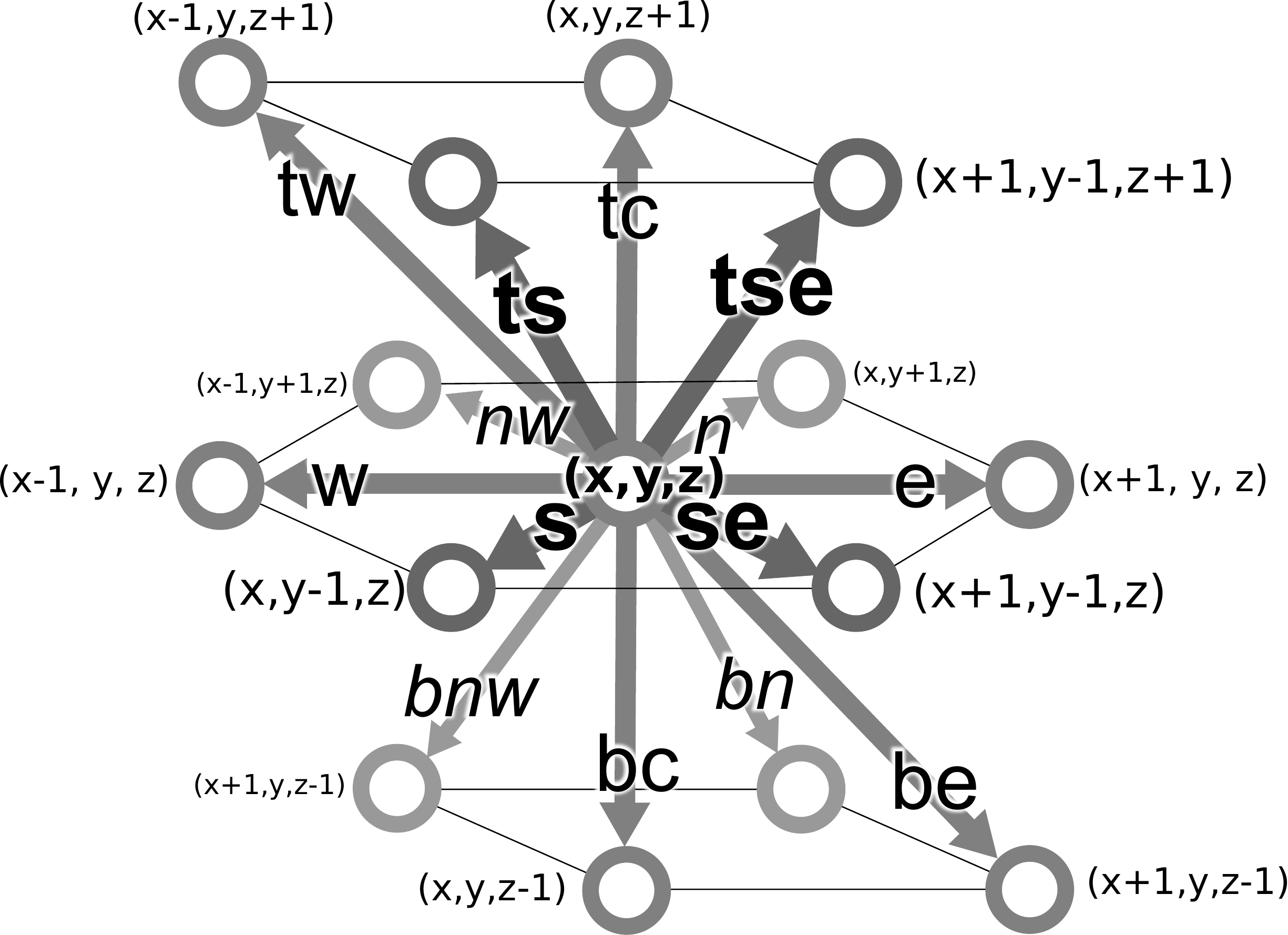}
\caption{
	Left: Direction vectors in a one-to-one correspondence between DoF and coordinates. 
	Right: Stencil directions and grid coordinates inside a structured tetrahedral grid.
	}
\label{fig:stencil-directions}
\end{figure}

Our surrogate ILU algorithm heavily relies on geometric properties associated with our DoF: 
Each micro-vertex in a macro-tetrahedron on level $L$ can be labeled by the logical grid coordinates 
\( G^{L}_{t} = \{ (x,y,z) \in \mathbb Z^3 \, : \, 0 \leq x, y, z \textmd{ and } x+y+z < 2^L+1\} \).
Similarly, we define the inner grid coordinates by \( \mathring{G}^{L}_{t} = \{ (x,y,z) \in \mathbb Z^3 \, : \, 1 \leq x, y, z \textmd{ and } x+y+z < 2^L\} \).
If we restrict the coordinates by setting z to a fixed value, we obtain a face-layer 
\( G^{N}_{f} = \{ (x,y) \in \mathbb Z^2 \, : \, 0 \leq x, y \textmd{ and } x+y < N \}. \)
For a vector $\bfu|_{I_t\cup\partial I_t}$ on level $l$ restricted to a tetrahedron $t$, there is a one-to-one correspondence between
DoF-indices in $\mcI_t $ and inner logical grid coordinates $\mathring{G}^{L}_t$ which can be constructed as follows:
Assume that $t$ is adjacent to the macro-vertices $v_i$ at coordinates $\tilde p_i \in \R^3$ for $1\leq i \leq 4$. 
The tetrahedron is spanned by the edges $d_i = \tilde p_{i+1}- \tilde p_1$ at the base point $\tilde p_1$ for $1\leq i \leq 3$ (see Figure~\ref{fig:stencil-directions} left).
The point $\tilde p^{(x,y,z)} = (d_1 \cdot x + d_2 \cdot y + d_3 \cdot z) / (2^L+1)$ for $(x,y,z) \in {G}^L_t$ belongs to a shape function $\phi_k \in V_{h_l}$ with $k \in \mcI_t \cup \mcI_{\partial t}$ such that
$\phi_k(\tilde p^{(x,y,z)}) = 1$. This induces the mapping $\iota_t: \mcI_t \cup \mcI_{\partial t} \to {G}^L_t$ with $\iota_t(k) = (x,y,z)$.
Thus, vector components $u_{i}$ of a vector $\bfu|_{\mcI_t\cup \mcI_{\partial t}}$ with $\iota \in \mcI_{t}\cup\mcI_{\partial t}$ will also be referred to by $u^{(x,y,z)}$ or $u^p$ for $p = (x,y,z) \in G^L_t$, when the macro-tetrahedron $t$ is evident from the context. 

The mapping between logical grid coordinates and local DoF in $\mcI_{t}$ also allows us to specify an ordering of the DoF indices. 
This is crucial since the properties of the Gauss--Seidel-Smoother (GS-Smoother) or the ILU-Smoother strongly depend on this order.
For $i, j \in \mcI_t$ with logical grid coordinates 
$(x_i,y_i,z_i) = \iota_t(i)$ and $(x_j,y_j,z_j) = \iota_t(j),$
we fix the ordering by 
$$i < j \implies (z_i < z_j) \vee (z_i = z_j \wedge y_i < y_j) \vee (z_i = z_j \wedge y_i = y_j \wedge x_i < x_j).$$
Consequently, the ordering strongly depends on the order of the adjacent macro-vertices $v_i$ which we used to construct $\iota_t$.
In Section~\ref{sec:3_hybrid_ilu_smoother}, we will use this by permutating $v_i$ with a permutation $\pi$ to obtain good smoothing factors $\mu_t$ for the new DoF ordering.

We will now introduce the applied stencil notation for our surrogate-ILU-Algorithm in Section~\ref{sec:4_matrix_free_ilu}.
For this, we first define the stencil directions between logical coordinates as displacement vectors, i.e. \( \{ x-y \,|\, x,y \in G^{L}_{t} \} \).
The most common directions are named after the four cardinal directions, as well as the top and bottom directions, such that the x-axis runs from west to east,
the y-axis from south to north, and the z-axis from top to bottom. For instance, the west direction $w$ corresponds to the displacement $(-1, 0, 0)$.
All stencil directions are collected in the set
\[ \mathcal D = \{ w, s, se, bnw, bn, bc, be, c, e, n, nw, tse, ts, tc, tw \}. \]
Relying on the ordering defined above, the set of all lower stencil directions needed in our ILU is given by 
\[ \mathcal D_l = \{ w, s, se, bnw, bn, bc, be \}. \]
Consider two indices $i \in \mcI_t$ and $j \in \mcI_t \cup \mcI_{\partial t}$ with coordinates $p_i = \iota_t(i)$ and $p_j = \iota_t(j)$.
Due to the local support of the low order conforming finite-element shape functions, we know that if $A_{ij} \neq 0$ there exists a $\tilde d \in \mcD$ such that $p_j = p_i + \tilde d$.
We can define the stencil $(A^{p_i}_{d})_{d \in \mcD}$ by $A^{p_i}_{\tilde d} = A_{ij}$.
The matrix-vector multiplication $v|_{\mcI_t} = (A u)|_{\mcI_t}$ on the macro-tetrahedron $t$ can therefore be written in terms of stencils as
\[
	v^p = \sum_{d \in \mcD} A^p_{d} u^{p+d}
	\quad \textmd{ for all } \quad p \in \mathring{G}^{L}_t,
\]
where we identified the DoFs with logical coordinates. 
Stencils and the associated grid coordinates are depicted in Figure~\ref{fig:stencil-directions}~(right).

We mainly rely on a geometric multigrid algorithm which combines a so called smoother with a coarse grid correction step to an optimal solver (see e.g.\cite{hackbusch2013multi}).
We now introduce smoothers which only act on the DoF of a single primitive. This is motivated by our hybrid mesh on which only operations between DoF located on the same primitive are cheap while everything else requires expensive inter-node communication. 
In our setting, for a given primitive $p \in \mcP_H$, a smoother acting on the DoF $\mcI_p$ located on the primitive can be described by applying a preconditioner matrix $C_p \in \R^{|\mcI_p|\times|\mcI_p|} $ inside a Richardson iteration with the appropriate restriction operators
\[
	\bfu \leftarrow \bfu + R^T_p C^{-1}_p R_p (\bff - A \bfu),
\]
where we denote the current right-hand-side vector by $\bff$ and the current estimate by $\bfu$.
Note that usually $C_p$ depends on our system matrix $A$ and is not necessarily symmetric.
If it is, we stress this by using $C_{p,sym}$.

As a reference example, we consider the GS-Smoother: On the primitive $p$ we define the lower triangular part $L_{A,p} \in \R^{|\mcI_p|\times|\mcI_p|}$ and diagonal part $D_{A,p} \in \R^{|\mcI_p|\times|\mcI_p|}$ of our system matrix $A$ restricted to the primitive as
\[
	(L_{A,p})_{ij} = \begin{cases}
		(R_p A R^T_p)_{ij} & \textmd{ if } j \leq i \\
		0 & \textmd{ else }
	\end{cases}
	\quad
	\textmd{ and }
	\quad
	(D_{A,p})_{ij} = \begin{cases}
		(R_p A R^T_p)_{ii} & \textmd{ if } i = j \\
		0 & \textmd{ else }
	\end{cases}
	.
\] 
The GS-Smoother is then given by $ C_{p} = L_{A,p} $ and its symmetrized version by 
$ C_{p,sym} = (L_{A,p}) (D_{A,p})^{-1} (L_{A,p})^T $.

\subsection{Steklov--Poincar\'e operator}\label{sec:2_2_steklov_poincare}

The idea of many non-overlapping domain decomposition methods is to solve a system of equations just on the boundary primitives of a decomposition and apply static condensation to the interior DoF.
Often, the decomposition is performed on unstructured meshes using graph partitioning libraries like METIS \cite{metis} or SCOTCH \cite{scotch}.
On hybrid meshes, the decomposition can simply be derived from the macro-mesh such that each refined macro-tetrahedron represents a subdomain. 
The macro-vertices, edges and faces and their respective DoF then form the interfaces between the subdomains.
The index set containing all DoF on the macro-interface is given by $\mcI_{\Gamma} = \cup_{t \in \mcT_H} \mcI_{\partial t}$.
With the previously introduced restriction operators, we define the submatrix $A_{\Gamma\Gamma} = R_{\Gamma} A R^T_{\Gamma}$ which just considers the coupling between the boundary DoFs, the matrix $A_{tt} = R_t A R^T_t$ which couples the interior DoF of a macro-tetrahedron $t \in \mcT_h$ and the coupling matrices
$A_{\Gamma t} = R_{\Gamma} A R^T_{t}$, $A_{t \Gamma} = R_{t} A R^T_{\Gamma}$ between the interior DoF of a tetrahedron and the boundary DoF.
Inverting $A$ is equivalent to solving $S \bfu_\Gamma = \boldsymbol{\chi}_\Gamma$ for $\bfu_\Gamma$ where the Steklov--Poincar\'e operator $S$ is given by
\[
    S =  A_{\Gamma\Gamma} - \sum_{t \in \mcT_H} A_{\Gamma{}t} A_{tt}^{-1} A_{t\Gamma}
\]
and the right hand side is  $\boldsymbol{\chi}_\Gamma = R_\Gamma \bfb - \sum_{t\in\mcT_H} A_{\Gamma t} A^{-1}_{tt} R_t \bfb$.
The interior DoFs can be reconstructed by solving $A_{tt} \bfu_{t} = R_t \bfb - A_{t\Gamma} \bfu_\Gamma$ for $\bfu_{t}$ on each tetrahedron $t \in \mcT_H$.
Usually, $S$ is not constructed explicitly, but is inverted by using a preconditioned conjugate gradient (PCG) method (see \cite{chan1994domain,toselli2004domain} for a non-exhaustive overview).
The PCG algorithms relies on the evaluation of matrix-vector products with $S$.
Evaluating $S$ requires the evaluation of $A_{tt}^{-1}$ for which we use a multigrid method with an ILU-Smoother in the interior of our macro-tetrahedron.
We will show, that our ILU algorithm gives us robustness with respect to tetrahedra which are distorted along one axial direction.

\section{Matrix-based ILU-Smoother}\label{sec:3_hybrid_ilu_smoother}
In this section, we introduce an ILU-\-Smooth\-er that is amenable to an efficient matrix-free implementation and investigate its performance.
Both the efficiency and the possibility for a matrix-free implementation heavily depend on the used ILU formulation.

We first start with the general definition:
Let the sparsity pattern of the matrix $A$ be given by
\( \mcS_A = \{ (i,j) \in \R^{|\mcI| \times |\mcI|} | A_{ij} \neq 0 \}, \)
then the ILU factorization consisting of the lower-triangular matrix $L\in \R^{|\mcI| \times |\mcI|} $ and diagonal matrix $D\in \R^{|\mcI| \times |\mcI|}$ is defined by 
\( (L D L^T)_{ij} = A_{ij} \) for $i,j\in \mcS_A$.

The ILU will be restricted to the interior of macro-tetrahedra.
Hence, we define the lower triangular matrices $L_t \in \R^{|\mcI_t|\times|\mcI_t|}$ and diagonal matrices $D_t \in \R^{|\mcI_t|\times|\mcI_t|}$ by
\( (L_t D_t L_t^T)_{ij} = (R_t A R^T_t)_{ij} \) for $i,j\in \mcS_{R_t A R^T_t}$.
We thus define our preconditioner on the macro-tetrahedra by $C_{t,sym} = L_t D_t L^T_t$.

\begin{remark}
\label{remark:modifiedILU}
Another popular choice would be the modified-ILU which is derived from the factorization $(\mathring{L} + D ) D^{-1} (\mathring{L} + D)^T_{ij} = A_{ij}$ for $(i,j) \in \mathcal{S}_{A}$, where $\mathring{L}$ is now a strictly lower-diagonal matrix. The remainder matrix containing the additional fill-in is $R = (\mathring{L} + D ) D^{-1} (\mathring{L} + D) - A$.
Summing over the remainder matrix terms and adding them to the diagonal with a weight $\omega$, i.e. $ (D_\omega)_{ii} = D_{ii} + \omega \sum_{i\neq j} \left| r_{ij} \right| $, yields the modified ILU
\[
    (\mathring{L} + D_{\omega}) D_{\omega}^{-1} (\mathring{L} + D_{\omega})^T,
\]
which for $\omega = 0$ becomes the usual ILU factorization, but is known to behave more robustly for different triangle types \cite{oertel1989multigrid}.
For $\omega > 0$, the $D_\omega$ matrix contains entries from the fill-in compared to the sparse matrix $A$, which we aim to approximate with the ILU.
In 3D, our fill-in consists of 12 additional nonzero entries per row which cannot be calculated in a memory-efficient way during our factorization.
They can possibly be reconstructed approximately in a postprocessing step, but we nevertheless restrict ourselves to the case of $\omega = 0$.
\end{remark}

For our concrete implementation, we aim for a factorization of the form $L D L^T$, where $L$ is a lower triangular matrix with a unit diagonal to minimize the number of multiplications and divisions during the forward and backward substitutions.

\subsection{Strategy on macro-elements}

The performance of the ILU smoother strongly depends on the ordering of the DoF.
To mitigate this effect and increase robustness, an alternating ILU-Smoother can be used \cite{oertel1989multigrid} which successively applies several ILU factorizations with different orderings.
This cannot be efficiently done in a matrix-free algorithm since the different orderings would result in cache misses for at least one of the orderings.
Thus, we have to use an efficient ordering from the beginning which means that we may have to permute the tetrahedral vertices in a preprocessing step.
It is shown in \cite{pinto2016robustness} with a Local Fourier Analysis (LFA) that for optimal performance, the triangles in 2D have to be orientated such that the first vertex is at the largest angle, the second at the smallest, and the third at the remaining angle.

This heuristic does not directly extend to 3D. 
For a more reliable strategy, we iterate over all macro-tetrahedra and apply an LFA with the techniques from \cite{gaspar2009fourier,gmeiner2013optimization} to the asymptotic ILU stencil for each of its possible orientations. 
Finally, we apply the permutation resulting in the smallest smoothing factor to the macro-tetrahedron.

For a formula of the asymptotic ILU stencil, we derive the ILU in-place factorization which only relies on local information from DoF neighbors.
On structured grids, the factorization can be derived by evaluating $L D L^T$ and comparing it to the matrix $A$ on the respective sparsity pattern. 
This yields the stencil equations
\[
  A^{p}_d
  = \sum_{\tilde d \in \mathcal D_l \cup \{c \}}
  L^{p}_{\tilde d} D^{p+\tilde d}_{c} L^{p+d}_{\tilde d - d}
  ,
  \quad
  \textmd{with}
  \quad
  d \in \mathcal D_l \cup \{c\}
  .
\]
If the stencils at $p+d$ for $d \in \mathcal D_l$ are already factorized, this system of equations can be solved for $L^{p}_{\tilde d}$ and $D^{p}_{c}$ for $\tilde d \in \mathcal D_l \cup \{ c \}$.
Thus, the factorization of the stencils $L^{p}_{\tilde d}$ and $D^{p}_c$ at $p = (x,y,z)$ at level $L$ relies on information at the points
\begin{align}
  I_p
  &
  =
  I^\beta_p \cup I^{\gamma}_p \quad\textmd{where}
  \\
  I^\beta_p &=
  \left\{
  (\tilde x, y, z) | 0 \leq \tilde x < x
  \right\}
  \cup
  \left\{
  (\tilde x, \tilde y, z) | 0 \leq \tilde x < 2^L\!+\!1\!-\!\tilde y\!-\!z, \tilde y < y
  \right\}
  \quad \textmd{ and }
  \\
  I^\gamma_p &=
  \left\{
  (\tilde x, \tilde y, z-1) |
  0 \leq \tilde x < 2^L\!+\!2\!-\!\tilde y\!-\!z,
  0 \leq \tilde y < 2^L\!+\!2\!-\!z
  \right\}
\end{align}
are the logical DoF-coordinates on the current layer and the layer below.
This standard procedure was already applied in the very first ILU paper\cite{meijerink1977iterative} for an incomplete Cholesky decomposition in 2D on quadrilateral grids.

Specifically, we obtain the following equations 
\begin{align}
\markblack{A^{p}_{bc}}
&=\markblack{L^{p}_{bc}} \markred{D^{p+bc}_{c}}
,
\\
\markblack{A^{p}_{s}}
&=\markblack{L^{p}_{bc}} \markred{D^{p+bc}_{c}} \markred{L^{p+s}_{bn}}
 + \markblack{L^{p}_{s}} \markred{D^{p+s}_{c}}
 ,
\\
\markblack{A^{p}_{bnw}}
&=\markblack{L^{p}_{bc}} \markred{D^{p+bc}_{c}} \markred{L^{p+bnw}_{se}}
 + \markblack{L^{p}_{bnw}} \markred{D^{p+bnw}_{c}}
 ,
\\
\markblack{A^{p}_{be}}
&=\markblack{L^{p}_{bc}} \markred{D^{p+bc}_{c}} \markred{L^{p+be}_{w}}
 + \markblack{L^{p}_{be}} \markred{D^{p+be}_{c}}
 ,
\\
\markblack{A^{p}_{w}}
&=\markblack{L^{p}_{bc}} \markred{D^{p+bc}_{c}} \markred{L^{p+w}_{be}}
 + \markblack{L^{p}_{bnw}} \markred{D^{p+bnw}_{c}} \markred{L^{p+w}_{bn}}
 + \markblack{L^{p}_{s}} \markred{D^{p+s}_{c}} \markred{L^{p+w}_{se}}
 + \markblack{L^{p}_{w}} \markred{D^{p+w}_{c}}
 ,
\\
\markblack{A^{p}_{bn}}
&=\markblack{L^{p}_{bc}} \markred{D^{p+bc}_{c}} \markred{L^{p+bn}_{s}}
 + \markblack{L^{p}_{be}} \markred{D^{p+be}_{c}} \markred{L^{p+bn}_{se}}
 + \markblack{L^{p}_{bnw}} \markred{D^{p+bnw}_{c}} \markred{L^{p+bn}_{w}}
 + \markblack{L^{p}_{bn}} \markred{D^{p+bn}_{c}}
 ,
\\
\markblack{A^{p}_{se}}
&=\markblack{L^{p}_{bc}} \markred{D^{p+bc}_{c}} \markred{L^{p+se}_{bnw}}
 + \markblack{L^{p}_{be}} \markred{D^{p+be}_{c}} \markred{L^{p+se}_{bn}}
 + \markblack{L^{p}_{se}} \markred{D^{p+se}_{c}}
 + \markblack{L^{p}_{s}} \markred{D^{p+s}_{c}} \markred{L^{p+se}_{w}}
\\
\intertext{and}
\markblack{A^{p}_{c}}
&=\markblack{D^{p}_{c}}
 + \left(\markblack{L^{p}_{bc}}\right)^{2} \markred{D^{p+bc}_{c}}
 + \left(\markblack{L^{p}_{be}}\right)^{2} \markred{D^{p+be}_{c}}
 + \left(\markblack{L^{p}_{bnw}}\right)^{2} \markred{D^{p+bnw}_{c}}
 + \left(\markblack{L^{p}_{bn}}\right)^{2} \markred{D^{p+bn}_{c}}
 \\ &\qquad
 + \left(\markblack{L^{p}_{se}}\right)^{2} \markred{D^{p+se}_{c}}
 + \left(\markblack{L^{p}_{s}}\right)^{2} \markred{D^{p+s}_{c}}
 + \left(\markblack{L^{p}_{w}}\right)^{2} \markred{D^{p+w}_{c}}
 ,
 \label{eq:matrix-equations}
\end{align}
at a grid point $p = (x,y,z)$,
where the already factorized symbols at points in $I_p$ are formatted in bold for convenience.
We further omitted the equations which require information that is not available during the factorization. 

To estimate the asymptotic stencils, we assume $K = \textmd{Id}_3$ inside the macro-tetrahedra.
The assumption is valid if $K$ is close to a scaled identity and it is not spatially varying too much.
Inside the macro-tetrahedra, this leads to coordinate independent stencils $A^p_d = A_d$ for all $d \in \mcD$ at all grid points $p$.
Similarly to the 2D case from \cite{pinto2016robustness}, we can now calculate the asymptotic ILU stencils $L^\infty_d$ and $D^\infty_c$ for $d\in\mcD$ by iterating Eqs.~\ref{eq:matrix-equations} with a Gauss-Seidel scheme.
The initial stencil values are given by $L^{(0)}_c = D^{(0)}_c = 1, L^{(0)}_d = D^{(0)}_d = 0$ for $d \in \mcD \setminus \{c\}$.

The Fourier symbols are then given by
\begin{align}
  D (\btheta) = D^{\infty}_c,  \quad
  L (\btheta) = 1 + \sum_{d \in \mcD_l} L^{\infty}_d e^{i d \cdot \btheta} \quad\textmd{ and } \quad
  A (\btheta) = \sum_{d \in \mcD} A_d e^{i d \cdot \btheta}
\end{align}
for $\btheta \in \left(-\pi, +\pi\right)^3$.
The smoothing factor $\mu_t$ is then
calculated for the highly-oscillating frequencies $\Theta_{osc} = \left(-\pi, +\pi\right)^3 \setminus \left(-\pi/2, +\pi/2\right)^3$ by
\[
  \mu_t = \sup_{\btheta \in \Theta_{osc}} \left|
    \frac{
    L(\btheta)
    D(\btheta)
    \conj{L}(\btheta)
    -
    A(\btheta)
    } {
    L(\btheta)
    D(\btheta)
    \conj{L}(\btheta)
    }
  \right|
  .
\]
For our purpose, $\mu_t$ is determined by sampling $\Theta_{osc}$ on a uniform grid with $16$ samples in each coordinate direction.

In order to determine the optimal orientation of a macro-tetrahedron $t$, we iterate over all possible vertex permutations $\pi_k$, $1 \leq k \leq 24$, 
calculate $\mu_{t_{\pi_k}}$ for all the permutated tetrahedra $t_{\pi_k}$ and finally apply the permutation with the smallest smoothing factor $\mu_{t_{\pi_{k}}}$.
Note that we need to calculate $16^3 \cdot 24 = 98304$ symbols per macro-tetrahedron to determine the final orientation, but we only have to do this once during a preprocessing step before the actual calculation.

\begin{remark}
We can a-priori determine the scaling behaviour of the asymptotic stencils in the mesh width $h$:
For this, assume that the asymptotic stencils $D^\infty_c$, $L^\infty_d$ for $d \in \mcD_l$ exist and that the stencil $A^p_d$ is independent of $p$.
Due to the transformation rule, we know $A_d \sim \Theta(h)$ and since $A_c = D^\infty_c ( 1 + \sum_{d \in \mcD_l}(L^\infty_d)^2) $ has just positive summands, we immediately see $D^\infty_c \sim \Theta(h)$,  $D^\infty_c (L^\infty_d)^2 \sim \Theta(h) $ and thus $L^\infty_d \sim \Theta(1)$.
\label{remark:asymptotic-relationships}
\end{remark}
Numerical results for the performance of our smoother for different permutations validating this strategy are given in Appendix~\ref{sec:appendix}.

\section{Matrix-free ILU-Smoother on macro-tetrahedra}\label{sec:4_matrix_free_ilu}
After introducing our matrix-based hybrid ILU-Smoother, we now turn towards our matrix-free algorithm.
Our goal is to develop a surrogate smoother with the same smoothing performance as the ILU, but with significantly lower memory requirements supporting large scale computations.
For this purpose, it is sufficient to consider structured meshes consisting of a single macro-tetrahedron since its behavior on hybrid grids is the same as for the matrix-based variant.
We also reflect this in our notation by omitting references to specific tetrahedra $t \in \mcT_H$ whenever possible and by keeping it as concise as possible.
Furthermore, we always assume that this macro-tetrahedron is oriented with respect to our reordering strategy.

\begin{figure*}[ht]
  \begin{center}
    \includegraphics[width=1.\textwidth]{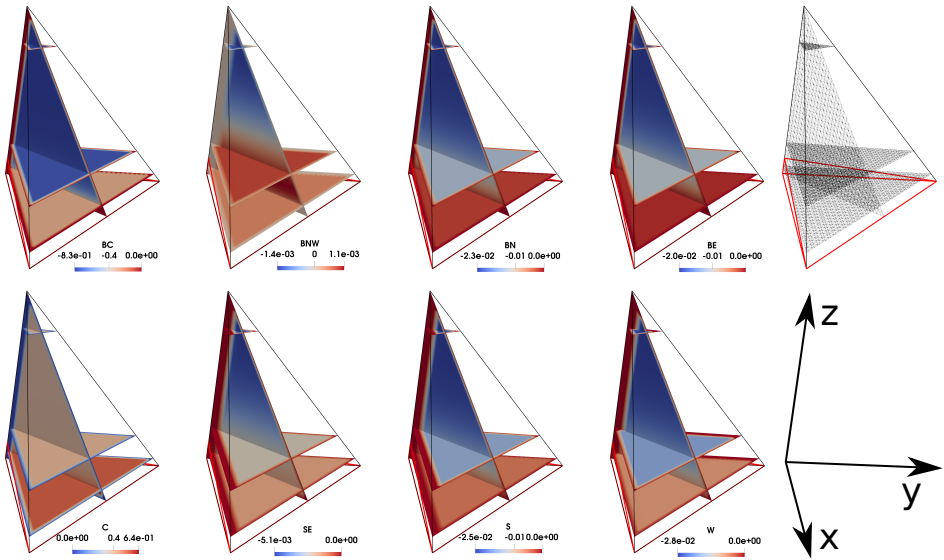}
  \end{center}
  \caption{\label{fig:surrogate-squished-stencils}
    Stencil functions for the $L^{(\cdot)}_d$ and $D^{(\cdot)}_c$ stencils of the distorted tetrahedron of height $h=0.1$ on grid level 5.
    The tetrahedron is scaled in the z-direction for better visibility. The original tetrahedron is depicted with red dashed lines in the plots.
  }
\end{figure*}
To motivate our approach, we depict in Figure~\ref{fig:surrogate-squished-stencils} each stencil direction of the factorized ILU of the distorted tetrahedron of height $h=0.1$ as continuous functions.
The z-axis is exaggerated by a factor of 10 to make the stencil plots easier to read.
The outlines of the distorted untransformed tetrahedron are indicated by red dashed lines in the plots while the outlines of the exaggerated tetrahedron are indicated by black lines.
One vertical and three horizontal slices make the inside visible.
For all stencil directions, we see that the stencil function does not vary much in the x-y-plane, except for some tiny layer close to the boundary.
Along the z-axis, all stencil functions display a color gradient from top to bottom.

This allows us to make several important observations:
Firstly, the stencil weights of a single direction are smooth functions. 
Secondly, the stencil weights are anisotropic with respect to the z-axis.
Both observations will enter into our surrogate strategy.

\begin{figure*}[h]
  \begin{center}
    \includegraphics[width=1.\textwidth]{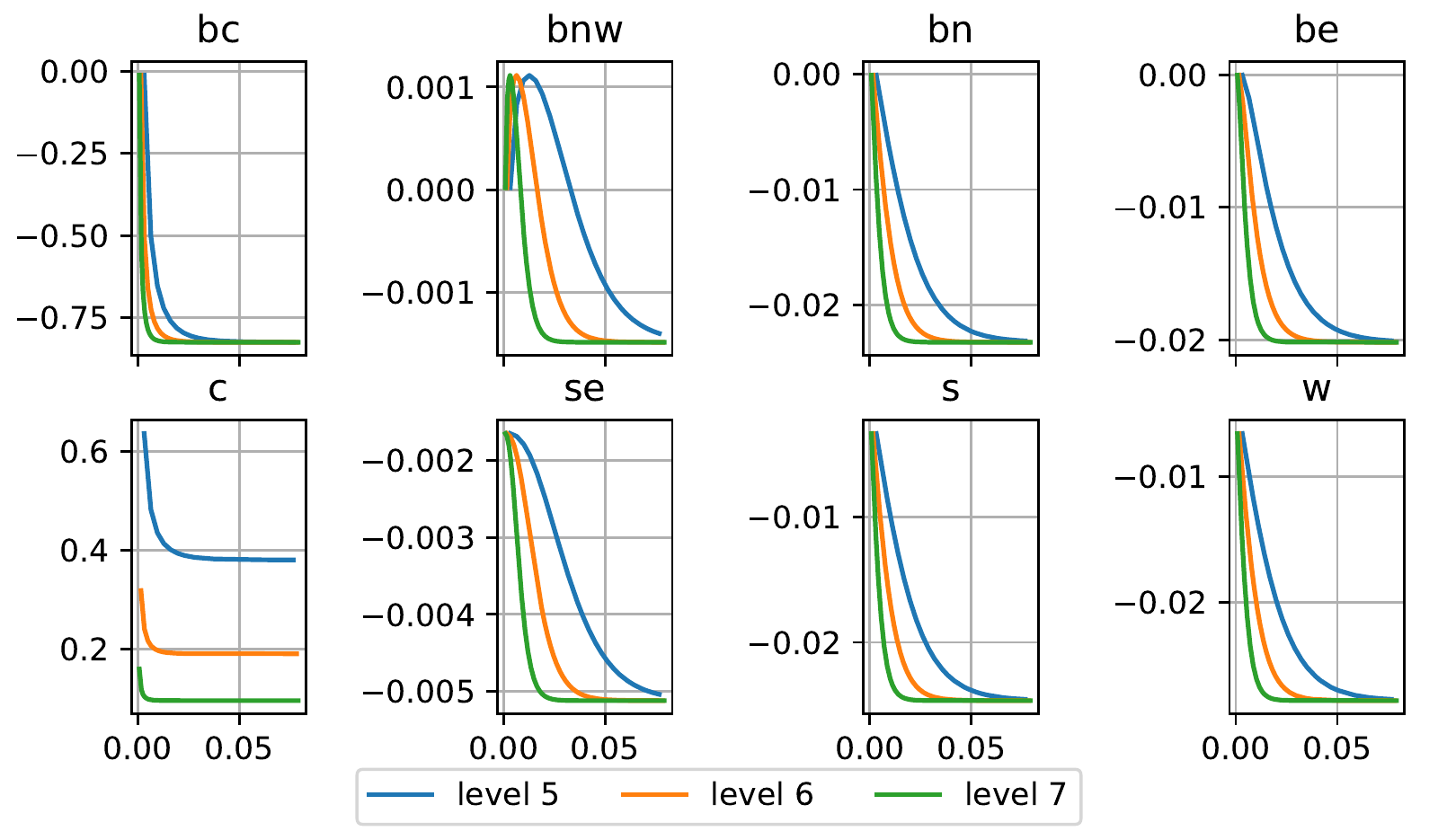}
  \end{center}
  \caption{\label{fig:surrogate-squished-stencils-plot}
    Stencil plots for the distorted tetrahedron on different grid levels along the line parallel to the z-axis starting at $(0.1,0.1,0)$ and ending at $(0.1, 0.1, 0.1)$.
  }
\end{figure*}
In Figure~\ref{fig:surrogate-squished-stencils-plot}, we have a more thorough look at the behavior in the z-direction.
We have depicted the stencil weights of the matrix-based ILU for the distorted tetrahedron along a line starting at $(0.1,0.1,0)$ and ending at $(0.1, 0.1, 0.1)$ for refinement levels from $5$ to $7$.
Each stencil direction approaches its asymptotic value towards the right.
From Remark~\ref{remark:asymptotic-relationships}, we expect that the diagonal stencil values from the $D$ matrix scale as the diagonal entries of the Laplacian with respect to the mesh width $h$.
Thus, its asymptotic value from a coarse level to its refinement has to decrease by a factor of $\nicefrac{1}{2}$ which is visible in the y-scaling of the center plots.
For the stencils of the $L$ matrix, Remark~\ref{remark:asymptotic-relationships} suggests that they do not scale in $h$ at all, which results in the same lower and upper bounds in the plots.
We see that all the curves for the different levels retain the same shape.
Furthermore, the curve for level 5 scaled by $\nicefrac{1}{2}$ or $\nicefrac{1}{4}$ in the x-direction yields curves similar to the ones on levels 6 and 7.
This suggests, that the non-asymptotic part of the stencil factorization decays in a level-independent way and only depends on the graph distance on the mesh to the DoF in the x-y-plane. 

\subsection{The surrogates}

We want to replace our in-memory ILU algorithm with a matrix-free variant to save main-memory storage capacities and enable large-scale computations.
Note that we cannot use an on-the-fly approach to execute the ILU Algorithm, since the factorization itself does not rely only on local information.
Furthermore, only the forward substitution moves through the DoF in the same direction as the typical in-place ILU factorization algorithm while the backward substitution moves in the opposite direction.
Rewriting the factorization such that it moves in the same direction as the backward substitution is possible, but it is not numerically stable due to the fact that it requires to reconstruct ILU stencils from nearly asymptotic ILU stencils.
Therefore, a memory-efficient factorization will only be the first step of our algorithm and has to be combined with some approximation approach that reconstructs the ILU.

For the factorization, we will use an in-place approach and store just enough information to complete the factorization.
Matrix entries which are no longer needed for the factorization are immediately discarded.
These types of in-place factorizations are typically used for ILU factorizations\cite{meijerink1977iterative} and just rely on local information from factorized neighboring DoF.
The factorization can be derived by solving Eqs~\ref{eq:matrix-equations} for $L^p_d, d \in \mcD_l$ and $D^p_c$ at each point $p = (x,y,z)$.
Our asymptotic storage requirements for the factorization can be estimated from the set $I_p$. 
Clearly, our memory requirements scale as $|I_p| = \mathcal{O}(h^{-2})$ for storing the factorization on the current face-layer ($I_\beta$) and the layer below ($I_\gamma$).
This has to be compared to a memory consumption scaling as $\mathcal{O}(h^{-3})$ which we would need to store the stencils on the whole grid in a matrix-based implementation.
The $\mathcal{O}(h^{-2})$ for the factorization are thereby negligible for practical purposes.

\begin{figure*}[ht]
  \begin{center}
    \includegraphics[width=0.75\textwidth]{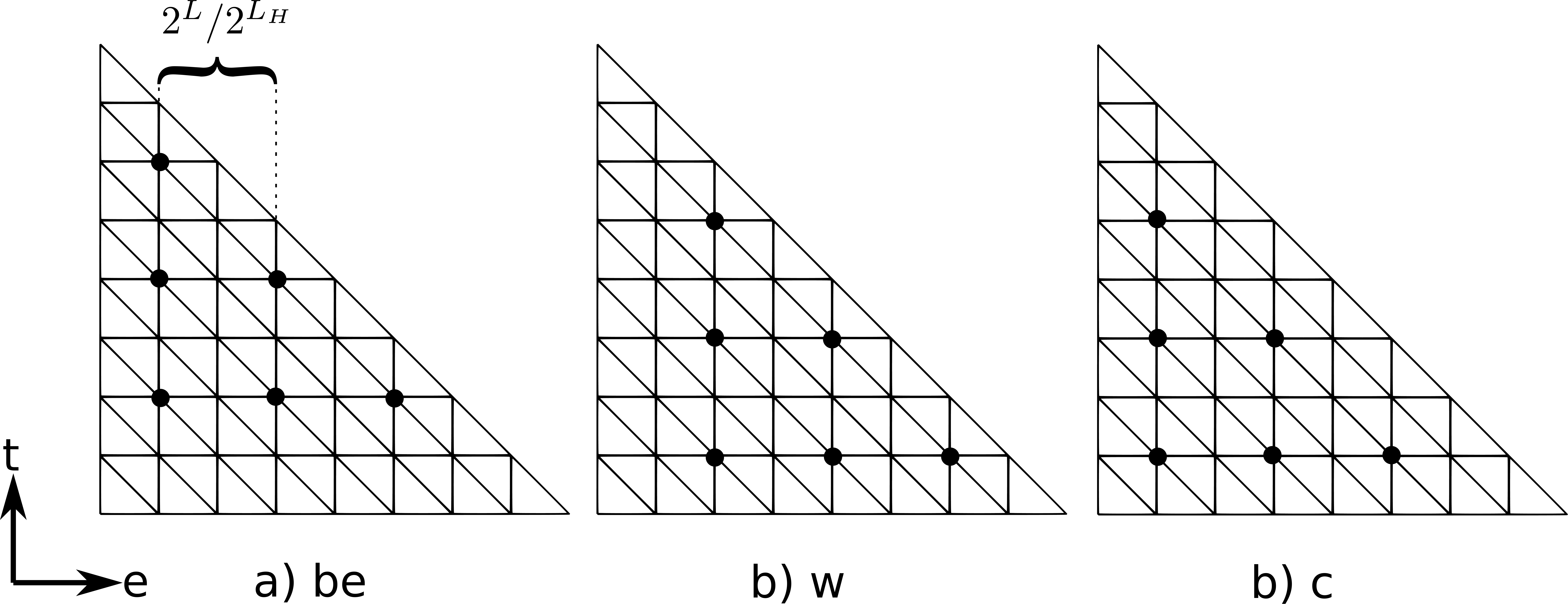}
  \end{center}
  \caption{\label{fig:sample-sets}
  Illustration of the sample sets in the x-z plane for $be$, $w$, and $c$.
  All points have the distance $2^L/2^{L_H}$ along each grid axis.
  Note the offset from the bottom for $be$ and the offset from the left for $w$.
  No shift has to be applied for the $c$ direction.
  }
\end{figure*}

We want to replace factorized stencils with surrogate polynomial approximations in the anisotropic polynomial space
\[
  \bbP_{\textmd{dg}_x, \textmd{dg}_y, \textmd{dg}_z}
  =
  \left\{
  x^i y^j z^k
  \,
  |
  \,
  0 \leq i \leq \textmd{dg}_x,
  0 \leq j \leq \textmd{dg}_y,
  0 \leq k \leq \textmd{dg}_z
  \right\}
  \,\textmd{for}\,
  \textmd{dg}_x,\textmd{dg}_y,\textmd{dg}_z \in \bbN
\]
such that $L_d^{(\cdot)}, (D^{-1}_c)^{(\cdot)} \in \bbP_{\textmd{dg}_x, \textmd{dg}_y, \textmd{dg}_z} $ for $d \in \mathcal D_l \cup \{ c \}$.
The different degrees for the different directions will allow us to vary the accuracy of the stencils in a direction-dependent way.
The factorized stencils on our grid will be added direction-wise to a least-squares (LSQ) problem to obtain one surrogate polynomial per stencil direction.
Of course, adding all the stencils to our least-squares problem would require a large memory overhead to store the matrix and a large computational overhead for solving the least-squares problem.
Therefore, only stencils $L^p_d$, $D^p_c$ at points $p \in S_d$ on some coarser grid $S_d$ of sample points will be added whose structure depends on the direction $d \in \mathcal D_l \cup \{ c \} $.
For a coarse sampling grid at level $L_H$, these sample points are defined as
\[
  S_d = \left\{
  p
  \,|\,
  p \in G^L_t 
  ,
  (p - (\mathbf{1}-d) )\,\textmd{mod}\, \min(2^L/2^{L_H}, 1)
  = 0
  \textmd{ and }
  p + d \not\in G^L_t \setminus \mathring G^L_t
  \right\}
  .
\]
They are illustrated for three directions in Figure~\ref{fig:sample-sets}.
The shift by direction $d$ is necessary to make sure to have the utmost boundary points in our approximation as sample points, which impacts the quality of the surrogate and thus the smoother.

\begin{algorithm}[htb]
  \begin{algorithmic}[1]
    \State $N = \textmd{ number of micro-vertices on a macro edge}$.
    \State Initialize $\beta^{(x,y)}_d = \gamma^{(x,y)}_d = 0$ for $d \in \mathcal D_l$ and $\beta^{(x,y)}_c = \gamma^{(x,y)}_c = 1$ and $(x,y) \in G_f^N$.
    \For{$z=1,\dots,N-2$}
    \For{$y=1,\dots,N-2-z$}
    \For{$x=1,\dots,N-2-z-y$}
    \State Assemble $A_d = A^{(x,y,z)}_d$ for $d \in \mathcal D$.
    \State Assign
    \State
    $ \begin{aligned}
        \markblack{\beta^{(x,y)}_{bc}}
         & =
        \markblack{A_{bc}}
        /
        \markred{\upgamma^{(x,y)}_{c}}
        \\
        \markblack{\beta^{(x,y)}_{s}}
         & =
        (
        \markblack{A_{s}}
        -\markblue{\beta^{(x,y)}_{bc}} \markred{\upgamma^{(x,y)}_{c}} \markred{\upbeta^{(x,y-1)}_{bn}}
        ) /
        \markred{\upbeta^{(x,y-1)}_{c}}
        \\
        \markblack{\beta^{(x,y)}_{bnw}}
         & =
        (\markblack{A_{bnw}}
        - \markblue{\beta^{(x,y)}_{bc}} \markred{\upgamma^{(x,y)}_{c}} \markred{\upgamma^{(x-1,y+1)}_{se}}
        ) / \markred{\upgamma^{(x-1,y+1)}_{c}}
        \\
        \markblack{\beta^{(x,y)}_{be}}
         & =
        (\markblack{A_{be}}
        -\markblue{\beta^{(x,y)}_{bc}} \markred{\upgamma^{(x,y)}_{c}} \markred{\upgamma^{(x+1,y)}_{w}})
        / \markred{\upgamma^{(x+1,y)}_{c}}
        \\
        \markblack{\beta^{(x,y)}_{w}}
         & =
        (
        \markblack{A_{w}}
        -\markblue{\beta^{(x,y)}_{bc}} \markred{\upgamma^{(x,y)}_{c}} \markred{\upbeta^{(x-1,y)}_{be}}
        - \markblue{\beta^{(x,y)}_{bnw}} \markred{\upgamma^{(x-1,y+1)}_{c}} \markred{\upbeta^{(x-1,y)}_{bn}}
        \\ & \qquad
        - \markblue{\beta^{(x,y)}_{s}} \markred{\upbeta^{(x,y-1)}_{c}} \markred{\upbeta^{(x-1,y)}_{se}}
        ) / \markred{\upbeta^{(x-1,y)}_{c}}
        \\
        \markblack{\beta^{(x,y)}_{bn}}
         & =
        (
        \markblack{A_{bn}}
        -\markblue{\beta^{(x,y)}_{bc}} \markred{\upgamma^{(x,y)}_{c}} \markred{\upgamma^{(x,y+1)}_{s}}
        - \markblue{\beta^{(x,y)}_{be}} \markred{\upgamma^{(x+1,y)}_{c}} \markred{\upgamma^{(x,y+1)}_{se}}
        \\ & \qquad
        - \markblue{\beta^{(x,y)}_{bnw}} \markred{\upgamma^{(x-1,y+1)}_{c}} \markred{\upgamma^{(x,y+1)}_{w}}
        ) / \markred{\upgamma^{(x,y+1)}_{c}}
        \\
                \markblack{\beta^{(x,y)}_{se}}
         & =
        (
        \markblack{A_{se}}
        -\markblue{\beta^{(x,y)}_{bc}} \markred{\upgamma^{(x,y)}_{c}} \markred{\upbeta^{(x+1,y-1)}_{bnw}}
        - \markblue{\beta^{(x,y)}_{be}} \markred{\upgamma^{(x+1,y)}_{c}} \markred{\upbeta^{(x+1,y-1)}_{bn}}
        \\ & \qquad
        - \markblue{\beta^{(x,y)}_{s}} \markred{\upbeta^{(x,y-1)}_{c}} \markred{\upbeta^{(x+1,y-1)}_{w}}
        ) / \markred{\upbeta^{(x+1,y-1)}_{c}}
        \\
        \markblack{\beta^{(x,y)}_{c}}
         & =
        \markblack{A_{c}}
        - \left(\markblue{\beta^{(x,y)}_{bc}}\right)^{2} \markred{\upgamma^{(x,y)}_{c}}
        - \left(\markblue{\beta^{(x,y)}_{be}}\right)^{2} \markred{\upgamma^{(x+1,y)}_{c}}
        \\ &\qquad
        - \left(\markblue{\beta^{(x,y)}_{bnw}}\right)^{2} \markred{\upgamma^{(x-1,y+1)}_{c}}
        - \left(\markblue{\beta^{(x,y)}_{bn}}\right)^{2} \markred{\upgamma^{(x,y+1)}_{c}}
        \\ &\qquad
        - \left(\markblue{\beta^{(x,y)}_{se}}\right)^{2} \markred{\upbeta^{(x+1,y-1)}_{c}}
        - \left(\markblue{\beta^{(x,y)}_{s}}\right)^{2} \markred{\upbeta^{(x,y-1)}_{c}}
        \\ &\qquad
        - \left(\markblue{\beta^{(x,y)}_{w}}\right)^{2} \markred{\upbeta^{(x-1,y)}_{c}}
      \end{aligned}$
      \State Set $p = h\cdot(x,y,z) \in \mathbb R^3$.
      \For{$d \in \mathcal D_l $}\label{alg:line:after-factorization}
      \State \textbf{if} $(x,y,z) \in S_d$ \textbf{then} Add $\beta^{(x,y)}_d$ at $p$ to a LSQ problem.
      \markteal{\small
      \State \textbf{if} $(x,y,z)$ on the cell boundary \textbf{then} Store $\beta^{(x,y)}_d$ in memory.\Comment{\bf(V1)} 
      }
      \EndFor
      \State \textbf{if} $(x,y,z) \in S_c$ \textbf{then} Add $1/\beta^{(x,y)}_c$ at $p$ to a LSQ problem.
      \markteal{\small
      \State \textbf{if} $(x,y,z)$ on the cell boundary \textbf{then} Store $1/\beta^{(x,y)}_c$ in memory.\Comment{\bf(V1)} 
      }
    \EndFor
    \EndFor
    \State Copy $\beta$ into $\gamma$.
    \State Set $\beta^{(x,y)}_c = 1$ and $\beta^{(x,y)}_d = 0$ for $d \in \mathcal D_l$ with $(x,y) \in G_f^{N-2-z}$.
    \EndFor
    \State Solve the least-squares problems.
  \end{algorithmic}
  \caption{
    \label{alg:3d-stencil-assembly} Memory efficient in-place $LDL^T$ factorization in a single tetrahedron.
    }
\end{algorithm}
This discontinuity also affects our substitution at the boundary.
For these stencils, we will test two variants:
In the first variant denoted by \markteal{\bf(V1)}, all the boundary stencils $L^p_d, D^p_d$ at boundary points $p$ will be stored and later used in the backward and forward substitutions.
Again, this just requires a memory complexity of $\mathcal{O}(h^{-2})$ since we have four boundary faces each with a $\mathcal{O}(h^{-2})$ storage requirement.
In a second variant \markviolet{\bf(V2)}, we will refrain from storing any boundary stencils and just correct the interpolated stencil entries at the boundary during the substitutions by setting the respective directions to 0.

The full factorization algorithm for 3D is given in Alg.~\ref{alg:3d-stencil-assembly}.
All the information in $I_p$ from previously factorized stencils is colored in red.
The information of the current stencil, which was calculated before and is used for other stencil entries, is colored in blue.
We denote temporarily saved stencils of the current face-layer belonging to $I^\beta$ by $\beta$ and stencils on the previous face-layer which correspond to $I^\gamma$ by $\gamma$.
After each factorization step, we check if the stencil direction should be added to a least-squares problem.
Note that for the central direction, we add the multiplicative inverse to the least-squares problem to avoid a floating-point division while applying the smoother.

Note that we rescale the integer coordinates $(x,y,z)$ by the mesh width $h$ before adding them to the least-squares problem in order to have a level independent approximation quality.
The steps in which we store the boundary stencils of the factorization are marked by \markteal{\bf(V1)} and they are only used in the first variant of the algorithm.
Before increasing the z index, we have to copy $\beta$ into $\gamma$.

\begin{algorithm}[h]
  \begin{algorithmic}[1]
    \For{$z = 1,\dots N-2$}
    \For{$y = 1,\dots N-2-z$}
    \For{$x =  1,\dots N-2-z-y$}
    \State Set $p = (x,y,z)$.
    \State Evaluate $A^{p}_d$ for $d \in \mathcal D$ by surrogates with Newton's Divided Differences
    \State$\quad$ Formula (NDDF) or assemble it.
    \State Apply $w^{p} = b^{p} - \sum_{d \in \mathcal D} A^{p}_d x^{p+d}$ for $d \in \mathcal D_l$.
    \If{$p$ is on the cell boundary}
    \State \markteal{ Load $L^{p}_d$ for $d \in \mathcal D_l$ from main memory. \Comment{\bf(V1)}}
    \State \markviolet{ If $p+d$ is on the boundary set $L^{p}_d = 0$. \Comment{\bf(V2)}}
    \Else
    \State Evaluate $L^{p}_d$ for $d \in \mathcal D_l$ by surrogates with NDDF.
    \EndIf
    \State $ w^{p} = w^{p} - \sum_{d\in \mathcal D_l} L^{p}_{d} w^{p+d} $
    \EndFor
    \EndFor
    \EndFor
    \For{$z = N-2,\dots 1$}
    \For{$y = N-2-z,\dots 1$}
    \For{$x =  N-2-z-y,\dots1, $}
    \State Set $p = (x,y,z)$
    \State Evaluate $(D^{p}_{c})^{-1}$ by surrogates with NDDF.
    \State $ w^{p} = w^{p} \cdot (D^{p}_{c})^{-1} $
    \For{$d \in \mathcal D_l $}
    \If{$p-d$ is on the cell boundary}
    \State \markteal{Load $L^{p-d}_{d}$ for $d \in \mathcal D_l$ from main memory. \Comment{\bf(V1)}}
    \State \markviolet{If $p-d$ is on the boundary set $L^{p-d}_{d} = 0$ for $d \in \mathcal D_l$. \Comment{\bf(V2)}}
    \Else 
    \State Evaluate $L^{p-d}_d$ for $d \in \mathcal D_l$ by surrogates with NDDF.
    \EndIf
    \EndFor
    \State $w^{p} = w^{p} - \sum_{d \in \mathcal D_l} L^{p-d}_{d} w^{p-d}$
    \State $x^{p} = x^{p} + w^{p}$
    \EndFor
    \EndFor
    \EndFor
  \end{algorithmic}
  \caption{
    \label{alg:3d-substitution-with-boundary}
    Surrogate evaluation for $x + (LDL^T)^{-1}(b - A x) $ with boundary layer.
  }
\end{algorithm}
In Alg.~\ref{alg:3d-substitution-with-boundary}, we show one step of the smoothing algorithm. 
Again, the steps specific to the first variant of the algorithm are marked by \markteal{\bf(V1)} while steps exclusively for the second variant are marked by \markviolet{\bf(V2)}.

Altogether, by merging the operator application with the forward substitution and the backward substitution with the diagonal scaling, we have to iterate twice over the entire mesh.
Therefore, the runtime due to the storage accesses should be comparable to the symmetric GS algorithm which also iterates twice over the tetrahedral grid.
During the first iteration, we calculate the residual by either assembling the $A$ matrix or evaluating its surrogate matrix.
We then plug the residual into the forward substitution for $L$.
In the second iteration, we combine the multiplication of the inverse diagonal matrix $D^{-1}$ with the backward substitution of $L^T$ and the correction step.
The evaluation of all the surrogate polynomials is implemented with a Newton's Divided Differences Formula (NDDF).
The algorithm is known to be unstable, however, this does not pose a problem for the low polynomial degrees considered here. This approach significantly speeds up the polynomial evaluation compared to other schemes: If the stencils for $L$ at $p = (x,y,z)$ were already evaluated, a new evaluation at the direct neighbor $(x+1, y, z)$ only takes 
$(\textmd{dg}_x + 1) \cdot |\mathcal D_l| = 7 \cdot (\textmd{dg}_x + 1)$ floating point additions and no multiplication at all.
Similarly, if during the second iteration the stencils for $L^T$ and $D^{-1}$ at $p = (x,y,z)$ were already evaluated, a new evaluation at the direct neighbor $(x+1, y, z)$ only takes 
$(\textmd{dg}_x + 1) \cdot (|\mathcal D_l| + 1) = 8 \cdot (\textmd{dg}_x + 1)$ floating point additions.
Since the number of stencils with a predecessor to the west is $\mathcal{O}(h^{-3})$ compared to $\mathcal{O}(h^{-2})$ stencils without a predecessor, the surrogate evaluation asymptotically becomes computationally cheap for small $h$ and large multigrid levels.
Furthermore, the stencils $L^p$, $D^p$ and $(L^T)^p$ taken together have the same number of non-zero stencil entries as $A^p$.
Therefore, we only have to evaluate one additional stencil compared to the symmetric GS-Smoother.
Note that the \emph{if} statement inside the loop keeps the pseudo-algorithm concise and in practice, it is avoided by moving it before the loop.
Therefore, for two sufficiently optimized implementations, we expect comparable computational costs for the ILU and the symmetric GS-Smoother.
In the following, we will show that we can preserve the advantageous convergence rates from Section~\ref{sec:3_hybrid_ilu_smoother} with our surrogate smoother.

\subsection{Numerical results}
In the following, all numerical experiments are conducted on a mesh hierarchy from levels 2 to 6.
We use a V-cycle based multigrid solver with 3 pre- and postsmoothing steps. A simple CG solver is used on the coarsest mesh. 
The asymptotic convergence factor $\rho$ is determined by 20-steps of a power iteration applied to the error propagation operator with a random initial error.

\begin{figure*}[ht]
  \begin{center}
    \includegraphics[width=1.\textwidth]{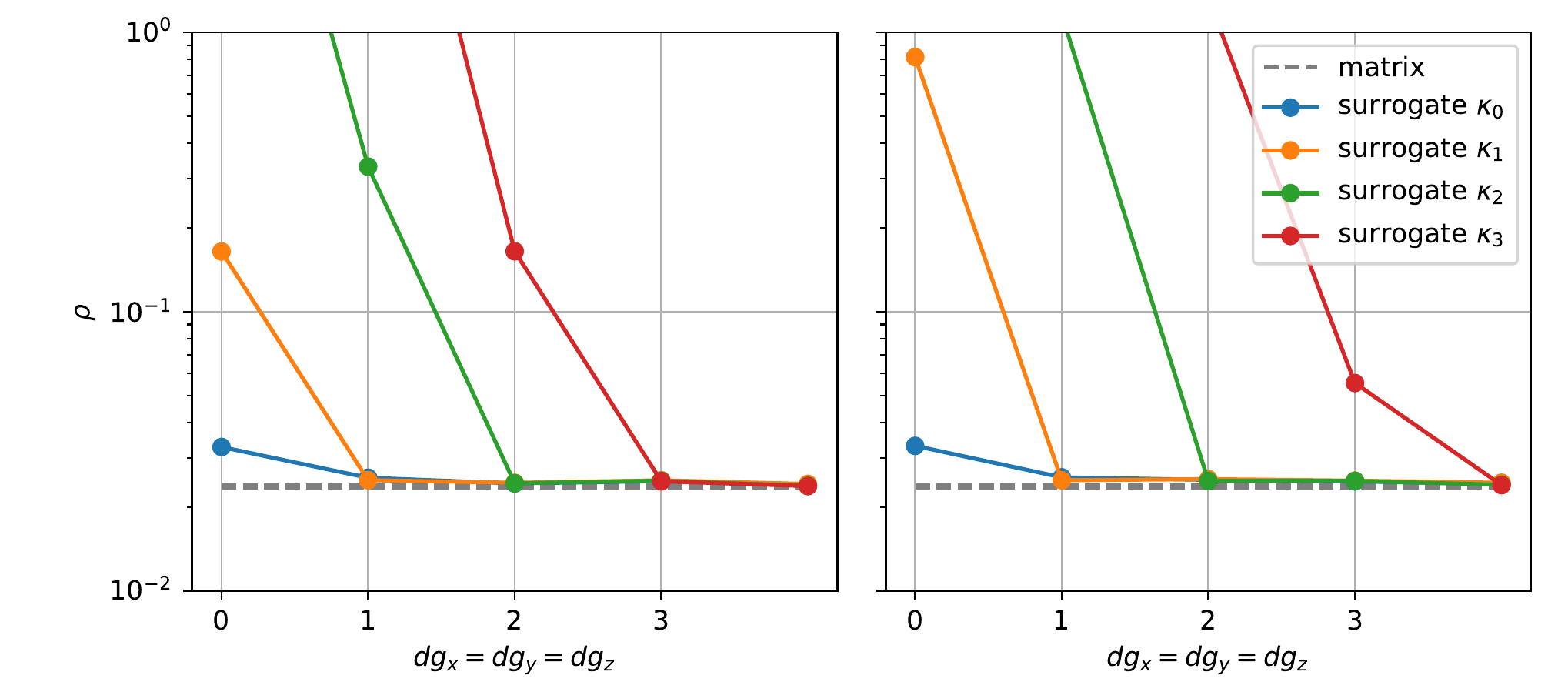}
  \end{center}
  \caption{\label{fig:surrogate-degrees-different-coefficients}
  Convergence rates of the surrogate ILU for different degrees for \markteal{\bf(V1)} on the left and \markviolet{\bf(V2)} on the right.
  }
\end{figure*}
To test the approximation quality of our ILU algorithm, we apply it to different scalar coefficients $\kappa_i$, where
\[
  \kappa_i (x, y, z) = 1 + 10 \left( x^i + y^i + z^i \right), \textmd{ with } 0 \leq i \leq 3
\]
is a polynomial of degree $i$ on a trirectangular tetrahedron with unit height and vertices $(0,0,0), (1,0,0), (0,1,0)$ and $(0,0,1)$.
The asymptotic convergence rates are given in Figure~\ref{fig:surrogate-degrees-different-coefficients} for \markteal{\bf(V1)} on the left and \markviolet{\bf(V2)} on the right.
The rates of the exact matrix version are plotted with dashed lines. 
The operator $A$ is approximated by a surrogate with $\textmd{dg}_x = \textmd{dg}_y = \textmd{dg}_z = 3$ while the degree of the surrogate approximation varies.
For undistorted tetrahedra, the convergence rates depend on the degree of the coefficient function $\kappa_i$.
Using the same degrees for the surrogate ILU and the coefficient function $\kappa_i$ allows us to recover the convergence rates of the original ILU algorithm.
Both versions of our algorithm perform similarly.
Saving the boundary stencils in \markteal{\bf(V1)} only provides a small advantage in case of $\kappa_3$ compared to using the surrogates everywhere and adjusting the boundary stencils in \markviolet{\bf(V2)}, see Figure~\ref{fig:surrogate-degrees-different-coefficients}.

\begin{figure*}[ht]
  \begin{center}
    \includegraphics[width=1.\textwidth]{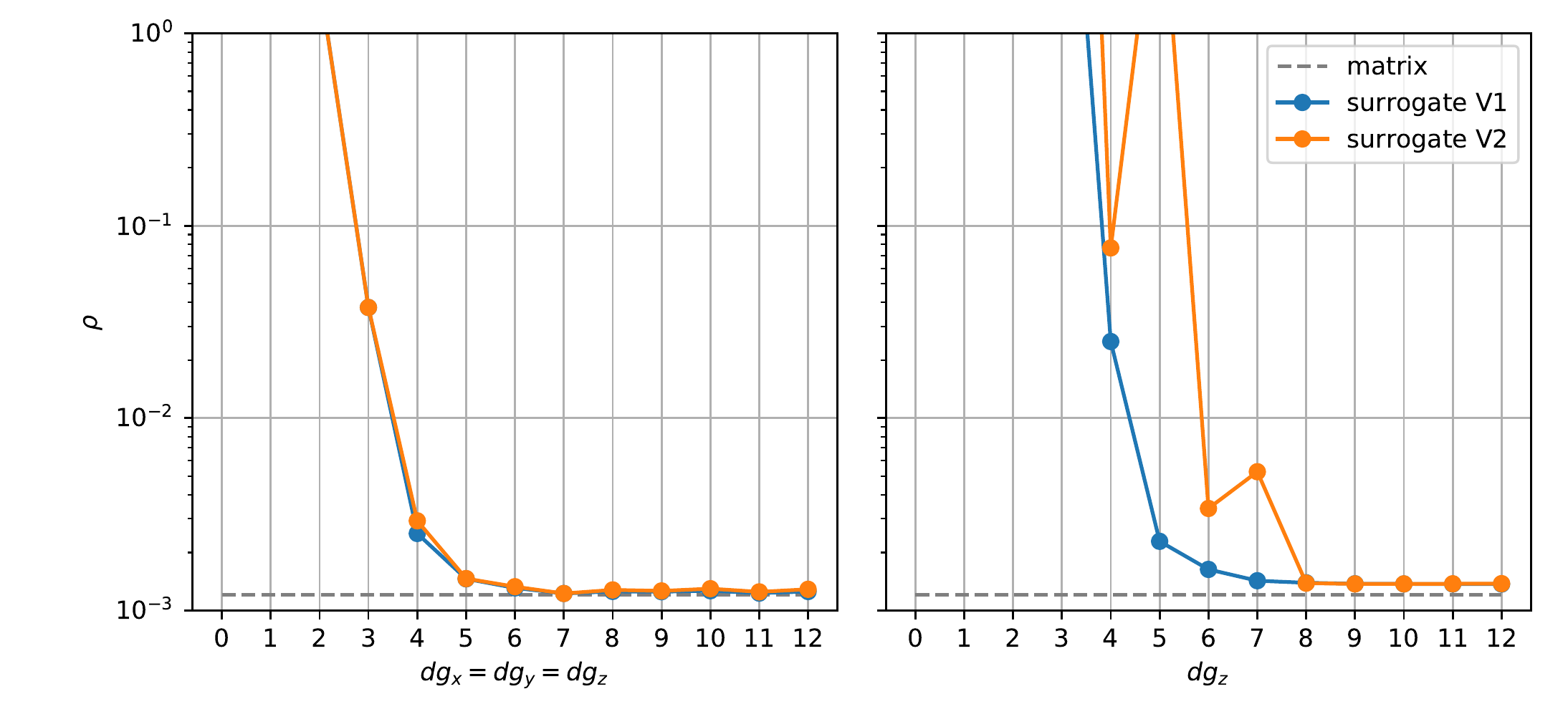}
  \end{center}
  \caption{\label{fig:surrogate-degrees-squished}
  Convergence rates of the surrogate ILU for different degrees for {\bf(V1)} on the left and {\bf(V2)} on the right.
  }
\end{figure*}
In Figure~\ref{fig:surrogate-degrees-squished}, we use the constant $\kappa = 1$ on a distorted tetrahedron of height $h=0.1$.
To approximate the Laplace operator for the residual, we use a constant surrogate polynomial. 
On the left, we use isotropic degrees along each coordinate axis for both variants of the surrogate ILU.
The convergence rates $\rho$ for both surrogate smoothers nearly coincide and saving the boundary stencils yields only negligible improvements.

The previous results from Figure~\ref{fig:surrogate-squished-stencils} suggest, that our stencil approximation has to be accurate in the z-direction, 
while lower polynomial degrees should be feasible inside the x-y-plane. 
Thus on the right of Figure~\ref{fig:surrogate-degrees-squished}, we set $\textmd{dg}_x = \textmd{dg}_y = 0$ and only vary $\textmd{dg}_z$.
In this case, during our iteration $\mathcal{O}(h^{-3})$ stencils do not have to be calculated at all, only the $\mathcal{O}(h^{-1})$ stencils along the z-axis have to be evaluated.
The asymptotic convergence rate of the correct ILU are not exactly attained even for large polynomial degrees, though the difference is negligible for all practical purposes.
To obtain a similar convergence rate as for the isotropic case, a higher polynomial degree in z-direction has to be used.
Saving the boundary stencils leads here to a more robust version, though again, both variants coincide if the polynomial degree in the z-direction is large enough.

\begin{figure*}[h]
  \begin{center}
    \includegraphics[width=1.\textwidth]{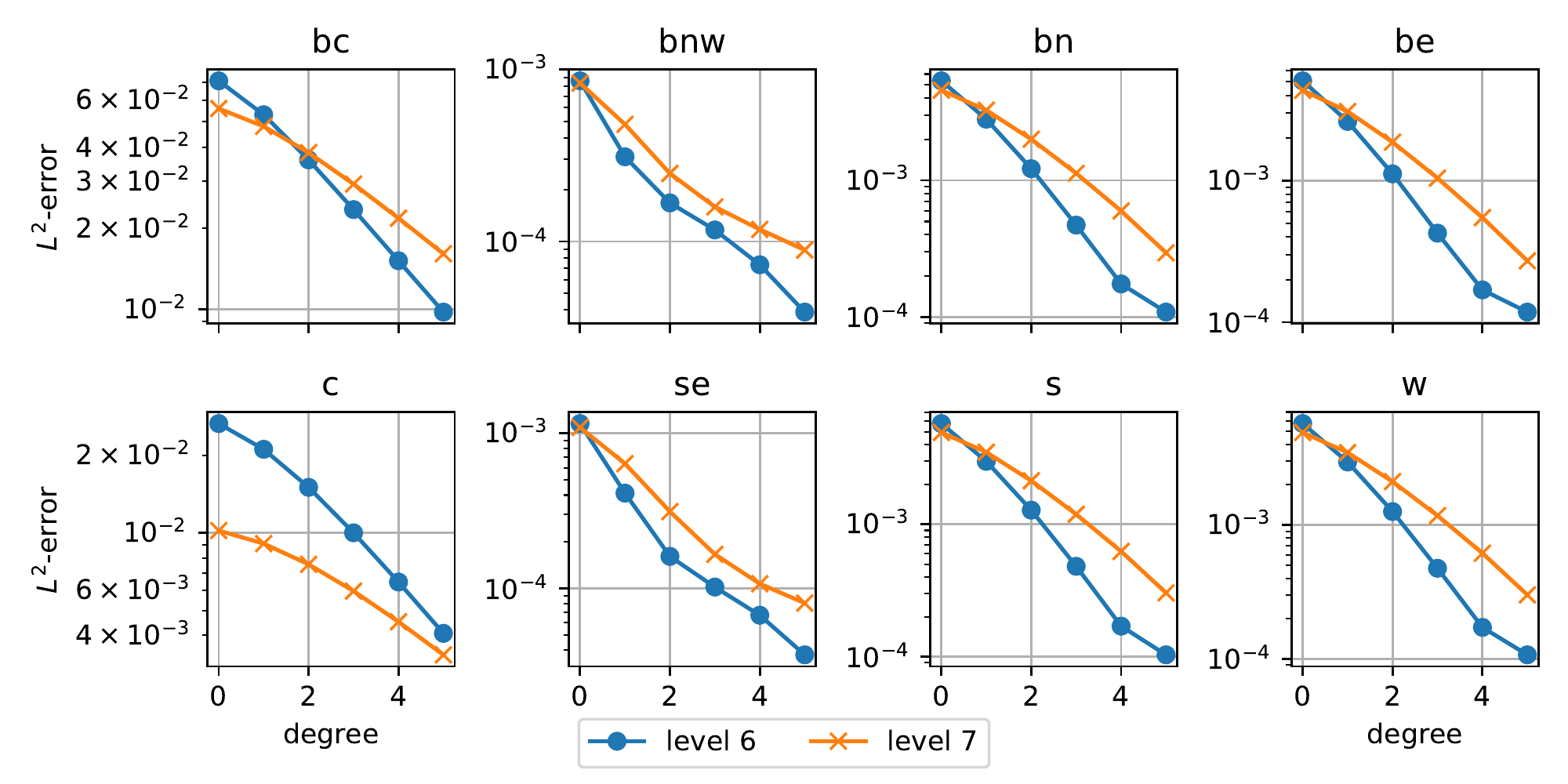}
  \end{center}
  \caption{\label{fig:surrogate-squished-stencils-L2}
    Error in the discrete $L^2$-norm for the $L^{(\cdot)}_d$ and $D^{(\cdot)}_c$ stencils of the distorted tetrahedron with height $h=0.1$ for levels $6$ and $7$.
    No values were skipped and no boundary stencils were saved.
  }
\end{figure*}
In Figure~\ref{fig:surrogate-squished-stencils-L2}, the error in the discrete $L^2$-norm is plotted on different multigrid levels for surrogate degrees from $0$ to $5$.
As expected, for all stencil directions, a higher polynomial degree leads to a smaller overall error.
For higher levels and hence finer meshes, this error becomes typically larger for the same polynomial degrees. 
Recall that in Figure~\ref{fig:surrogate-squished-stencils-plot}, we observed that the non-asymptotic parts of the factorization are scaled in the spatial direction by a factor of $\nicefrac{1}{2}$ from the coarser to the finer level. 
Thus, the stencil functions which we aim to approximate obtain a larger gradient and hence require higher polynomial degrees for the same approximation error.
An exception is the c-direction, whose absolute value (Figure~\ref{fig:surrogate-squished-stencils-plot}) decreases by a factor of $\nicefrac{1}{2}$ from the coarser level to the finer level. 
This results in a smaller absolute error. 

\begin{figure*}[ht]
  \begin{center}
    \includegraphics[width=.49\textwidth]{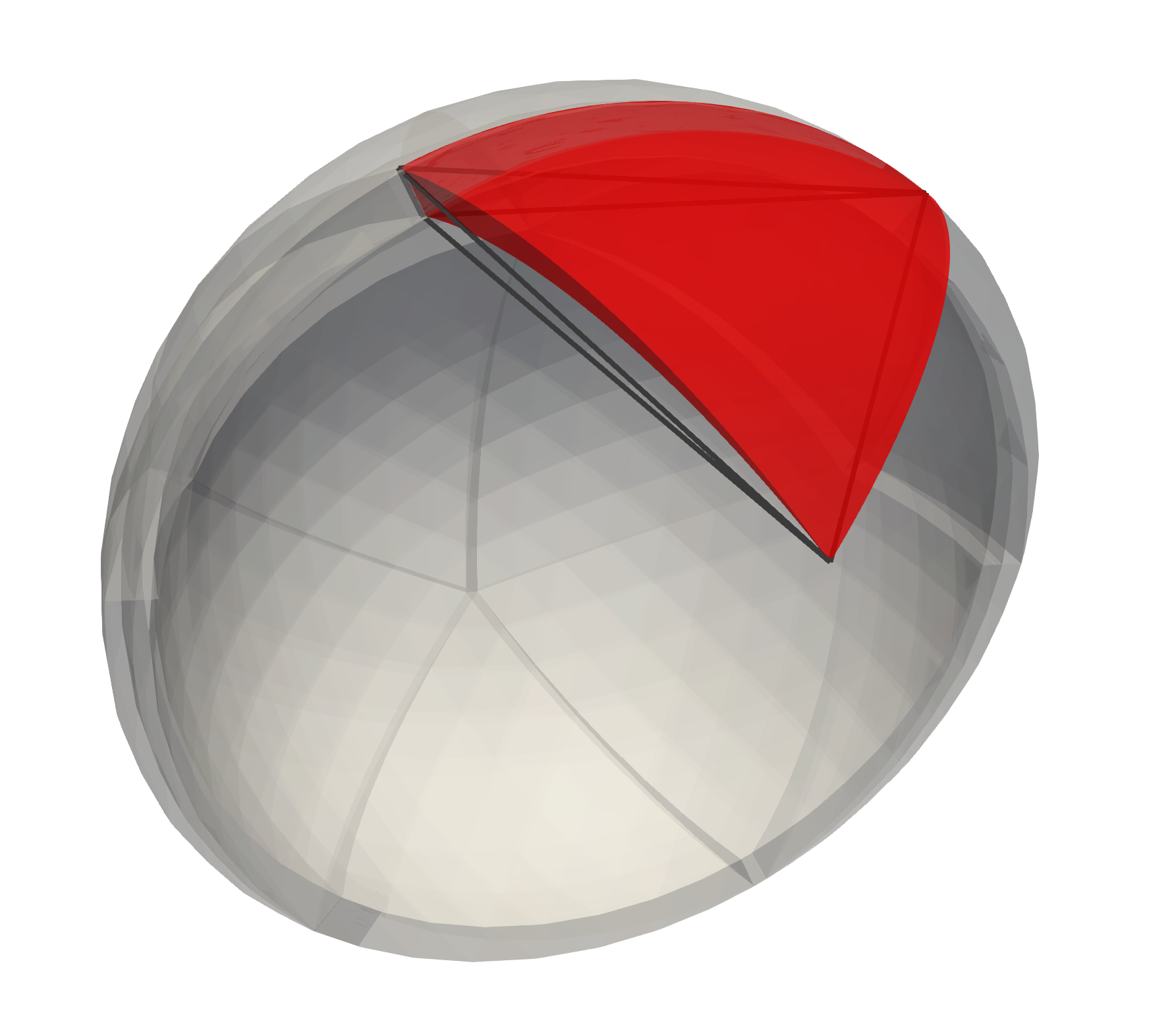}    \includegraphics[width=.49\textwidth]{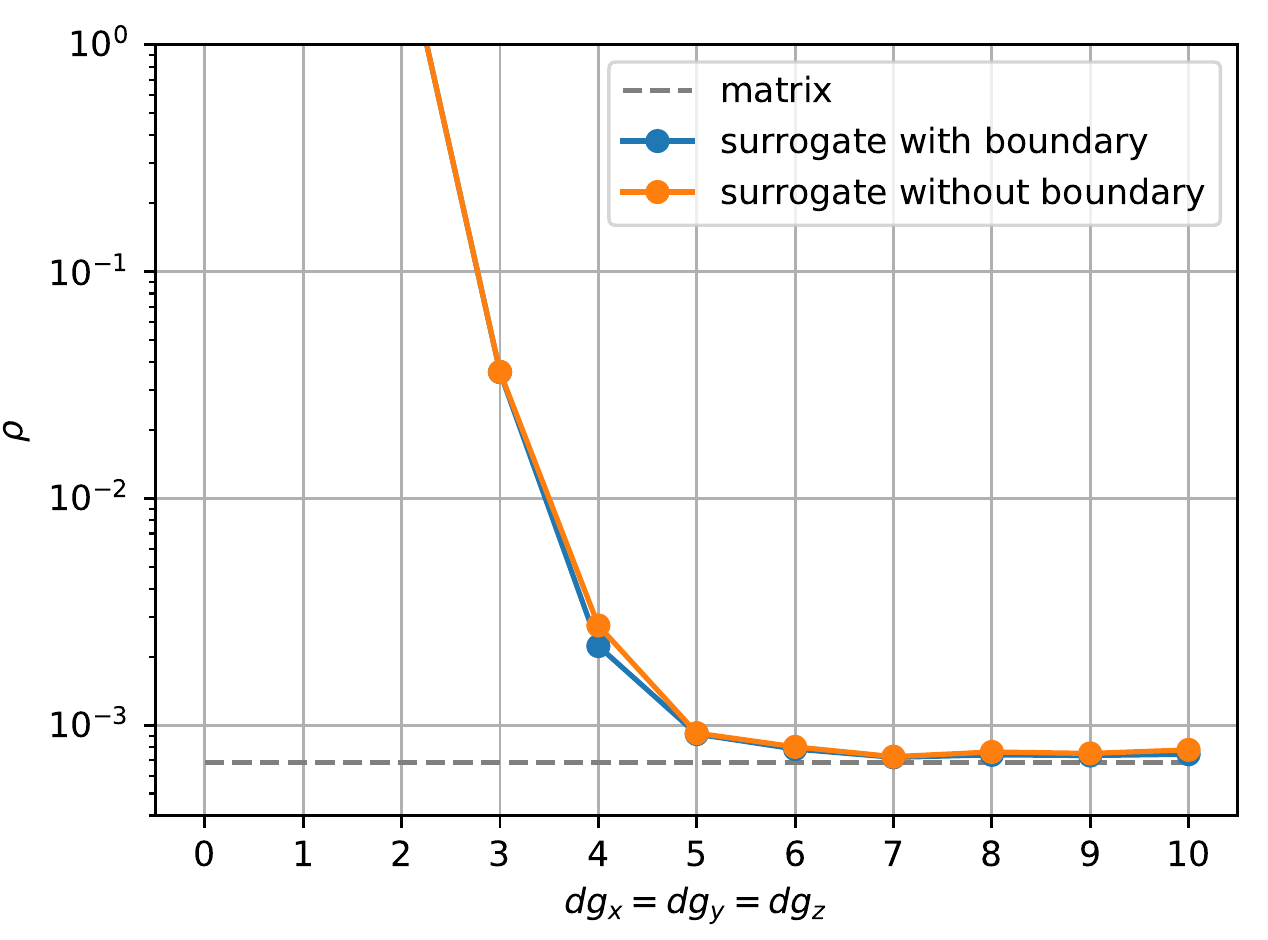}
  \end{center}
  \caption{\label{fig:surrogate-degrees-blending}
  Left: Test geometry with a distorted and blended tetrahedron.
  Right: Convergence rates of the surrogate ILU on the given primitive for different degrees.
  }
\end{figure*}
To show that our surrogate ILU also works in scenarios with less artificial geometries and coefficients, we apply it to a blended tetrahedron on the outer boundary of a narrow spherical shell with outer radius 1 and inner radius 0.9 (Figure~\ref{fig:surrogate-degrees-blending} left).
The blended tetrahedron under consideration is depicted in red while the unperturbed tetrahedron is indicated with black lines.
The rest of the shell which we do not use in the simulation is indicated in light-gray.
This benchmark combines both properties of the previous example: the varying coefficient tensor due to the blending map and the strong anisotropy due to the distorted base triangle. 

For simplicity, we reassemble $A$ for each operator application instead of relying on a surrogate approximation.

The convergence rates are depicted on the right of Figure~\ref{fig:surrogate-degrees-blending}, for surrogate polynomials with isotropic degrees along the coordinate axes.
Starting with polynomial degree~$3$, our ILU-Smoother is accurate enough to ensure convergence.
In this case, both variants of our ILU-Algorithm behave in the same way.
For polynomials of degrees larger than $7$, the rates for the surrogate ILU and the matrix-based ILU coincide.

\begin{figure}
    \centering
    \includegraphics[width=0.75\textwidth]{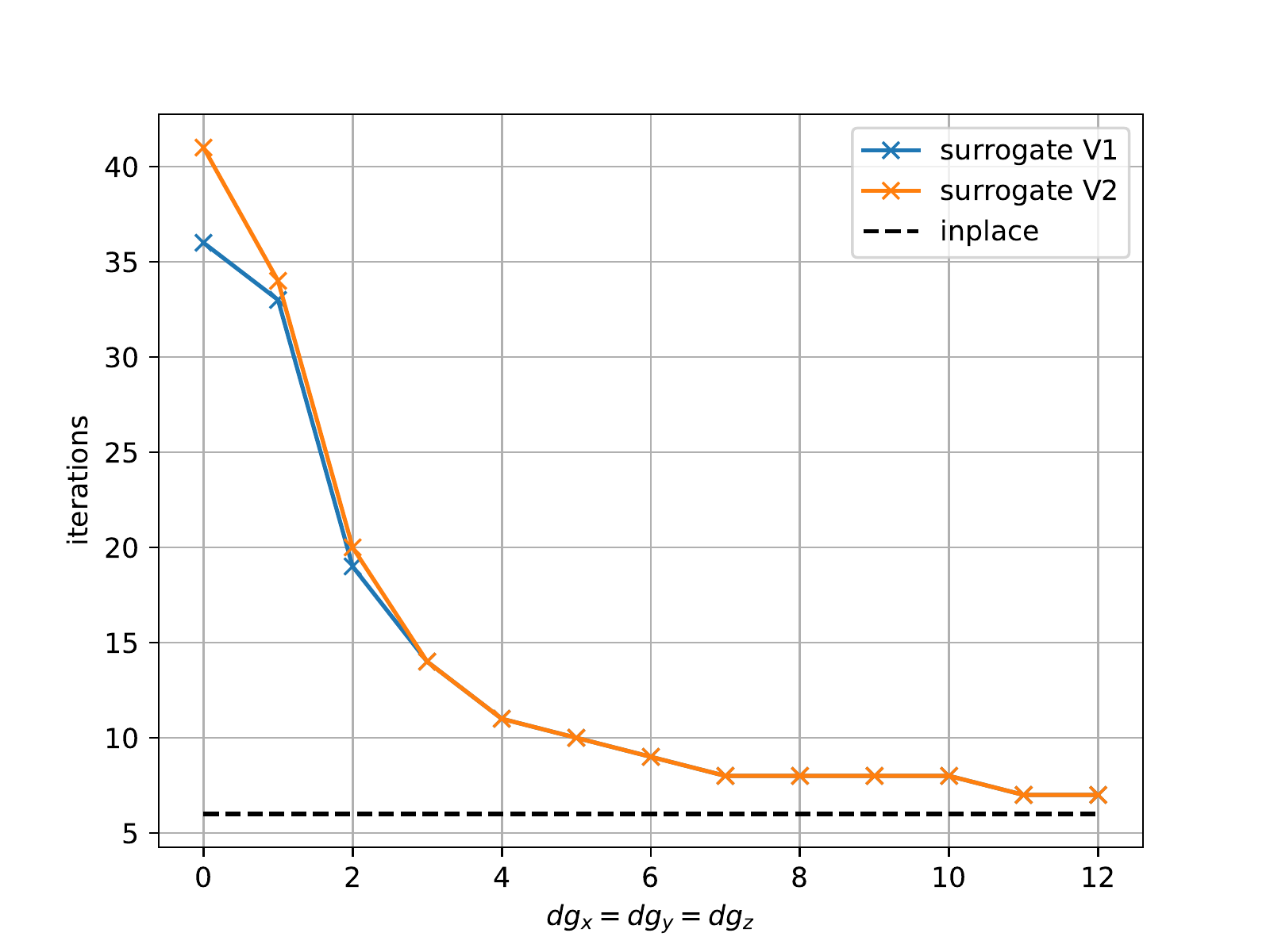}
    \caption{
      Number of required iterations for a CG method directly preconditioned by our ILU factorization to decrease the absolute residual below $10^{-3}$ for the example from Fig.~\ref{fig:surrogate-degrees-blending}.
    }
    \label{fig:ilu-preconditioned}
\end{figure}

An ILU-factorization is not only used as a smoother for multigrid methods, but can also be applied directly as a preconditioner for a different iterative method. 
Candidates would be the Richardson iteration or Krylov methods.
This use-case is depicted in Figure~\ref{fig:ilu-preconditioned}: 
We apply the surrogate ILU as a preconditioner within a CG method and increase its degree along the x-axis. 
We present the required number of iterations for decreasing the absolute residual below $10^{-3}$.
The symmetric Gauss--Seidel method, not depicted in the plot, needs $67$ iterations for achieving this.
The ILU's performance is much better, starting at $41$ and $37$ iterations for low polynomial degrees and resulting in $9$ to $7$ iterations for higher degrees.
Similar to the multigrid case, a surrogate polynomial degree of $7$ is roughly the point at which our approximation quality stagnates and does not improve by much. 

\subsection{Performance analysis}

\begin{table}[ht]
    \centering
    \begin{tabular}{c|cc|cc}
    & \multicolumn{2}{c|}{matrix-free} & \multicolumn{2}{c}{in-memory}
    \\ \hline
    line & memory traffic {\small[Byte]} & FLOP & memory traffic {\small[Byte]} & FLOP \\
    \hline
    5 & 0 & 0 & 0 & 0 \\
    6 & $15 \cdot 8$ & 16 & $15 \cdot 8$ & 16  \\
    11 & 0 & $7 \cdot \textmd{deg}_x$ & $7 \cdot 8$ & 0\\
    12 & $6 \cdot 8$ & 8 & $6 \cdot 8$ & 8\\
    \hline
    $\sum$ & $21 \cdot 8$ & $24 + 7 \cdot \textmd{deg}_x $
    & $28 \cdot 8$ & 24 
    \\
    \hline
    \multicolumn{1}{l|}{AI} 
    & \multicolumn{2}{r|}{ $(\nicefrac{1}{7} + \nicefrac{\textmd{deg}_x}{24})$ FLOP/Byte } 
    & \multicolumn{2}{r}{ $\nicefrac{3}{28}$ FLOP/Byte} 
    \\
    \hline
    17 & 0 & $1 \cdot \textmd{deg}_x$ & $1 \cdot 8$ & 0 \\
    18 & $1 \cdot 8$ & 1 & $1 \cdot 8$ & 1 \\
    24 & 0 & $7 \cdot \textmd{deg}_x$ & $7 \cdot 8$ & 0 \\
    25 & $6 \cdot 8$ & 8 & $6 \cdot 8$ & 8 \\
    26 & $1 \cdot 8$ & 1 & $1 \cdot 8$ & 1 \\
    \hline
    $\sum$ & $8 \cdot 8$ & $10 + 8 \cdot \textmd{deg}_x$ & $16 \cdot 8$ & 10 \\
    \hline
    \multicolumn{1}{l|}{AI} 
    & \multicolumn{2}{r|}{ $(\nicefrac{5}{32} + \nicefrac{\textmd{deg}_x}{8})$ FLOP/Byte}
    & \multicolumn{2}{r}{ $\nicefrac{5}{64} $ FLOP/Byte} \\
    \end{tabular}
    \caption{Theoretical cost analysis of memory traffic and FLOP in the interior of our mesh.}
    \label{tab:cost_surrogate}
\end{table}
Large scale low-order finite-element computations are typically memory-bound.
Therefore, to increase the performance of an algorithm we have to minimize memory transfers.
To judge the possible performance benefits of our surrogate algorithm, we will therefore determine the arithmetic intensity of Algorithm~\ref{alg:3d-substitution-with-boundary} away from the boundary.
The costs for the different lines of our implementation are given in Table~\ref{tab:cost_surrogate}. 
The columns on the left relate to our matrix-free implementation while the right column lists the theoretical costs for a matrix-based implementation.
We separately sum up the costs of the first and second loop in Algorithm~\ref{alg:3d-substitution-with-boundary}.

In this analysis, we assume that no blending and a constant coefficient is used. Hence, the matrix stencil $A^p_d$ is constant, fits into the caches and therefore Line~6 does not require memory traffic or FLOP.
In Line~7, we can assume that $x^{p+w}$ is already present in the cache.
Therefore, only $15$ vector entries of $x$ and $b$ have to be loaded, and $15$ multiplications and one subtraction have to be applied.
Since we assume to be in the interior, Lines~7-10 do not contribute.
In Line~11, evaluating each of the 7 stencil entries takes $\textmd{deg}_x$ evaluations.
In Line~12, we assume that $w^{p}$ and $w^{p+w}$ are already cached and that 7 multiplications followed by 1 addition are needed.

The in-memory implementation differs from the matrix-free implementation in Lines 11 and 24, where instead of evaluating the stencil matrix-free, it has to loaded from main memory. Altogether, $15 \cdot 8$ additional bytes have to be loaded from the main memory. Therefore, our matrix-free implementation has the potential to decrease the necessary memory traffic by $34 \%$.

Finally, we want to discuss the total memory requirements of the factorization.
Both for the matrix-free and the in-memory algorithm, we need one additional intermediate vector $w$ which consists of $\textmd{\#dof} \cdot 8 $ bytes.
Since the memory requirements for the polynomial coefficients can be neglected, no additional memory for the matrix free implementation is required.
The in-memory algorithm on the other hand requires the storage of the matrices $L$ and $D$. 
Since each interior row of $L$ consists of 7 non-zeros, in total $\textmd{\#dof} \cdot 64 $ additional bytes have to be stored.
Therefore, the in-memory algorithm requires 9 times the memory of our matrix-free implementation. 

\begin{figure}
    \centering
    \includegraphics[width=0.75\textwidth]{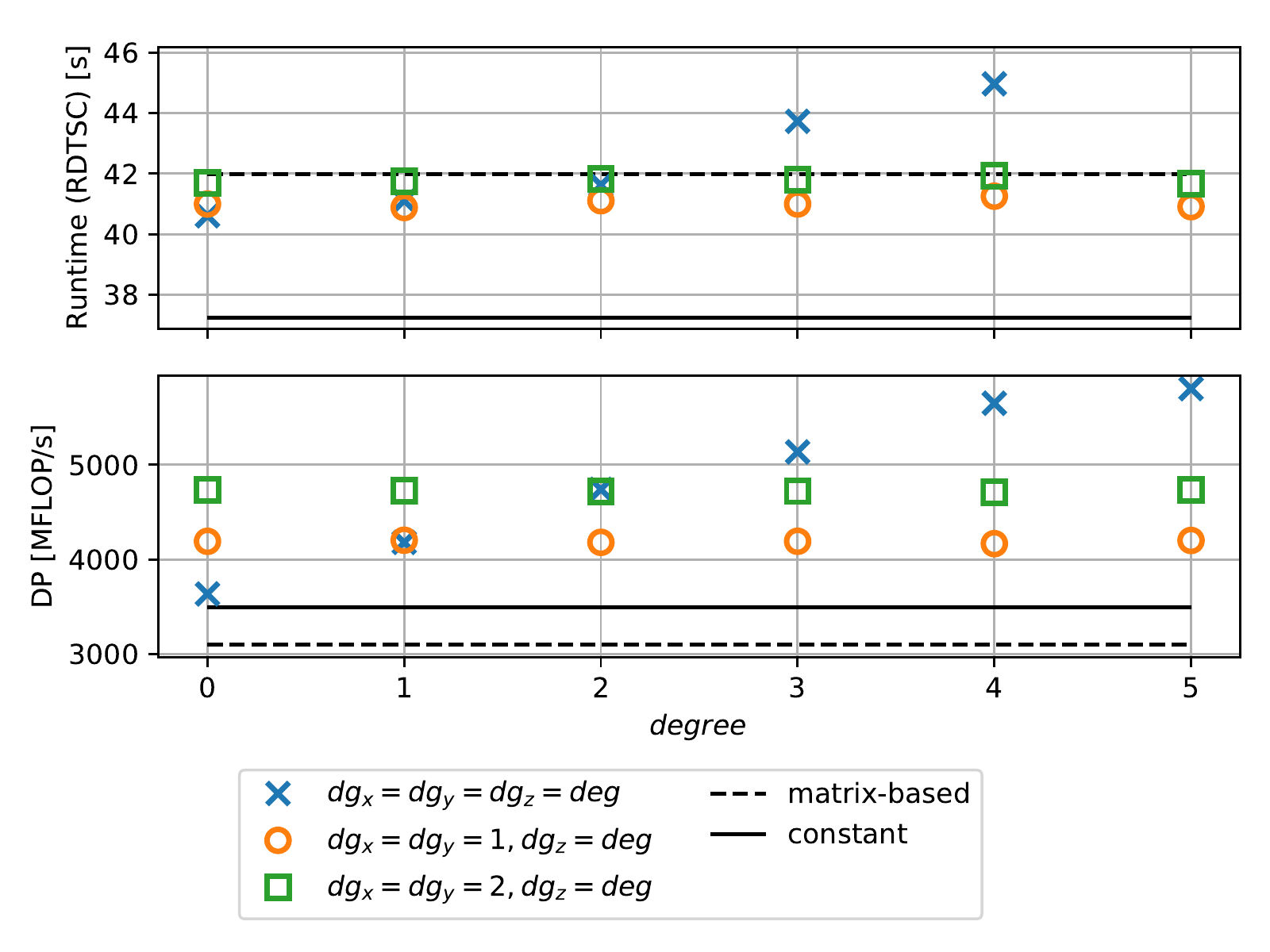}
    \caption{
    Performance for the forward-substitution:
    Top: Runtime for different surrogate polynomials.
    Bottom: Scalar double-precision FLOP/s for different surrogate polynomials.
    }
    \label{fig:performance}
\end{figure}
In Figure~\ref{fig:performance}, we depict the experimentally determined runtime and the FLOP/s averaged over the processes of the second variant of our forward substitution.
The code is executed on an Intel Xeon Gold 6136 with 12 processes distributed over its two sockets.
As an example, we used a grid with 48 macro-tetrahedra on multigrid level 9 with 32 repetitions of the forward- and backward substitutions,
and as a measurement tool we rely on LIKWID \cite{7103452,psti}. 
In the figure, the performance of a simple matrix-based implementation, in which all the stencils of the $L$, $D$ and $L^T$ matrices are stored continuously inside a large memory buffer is depicted with a dashed line while the results of an implementation which just relies on fixed asymptotic stencils is denoted with a solid line.
The latter provides a lower bound for the attainable runtime of our algorithm, but has no practical relevance.

If we increase the polynomial degree in each direction uniformly, we arrive at the {\includegraphics[width=6pt]{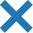}} curve.
With each increase in the degree, the FLOP count increases, and brings a modest increase in runtime.
For a degree larger or equal 3, storing the stencils externally and loading them from memory is more efficient.

We showed that usually a high degree in the z-direction and lower degrees in the x- and y-directions are completely sufficient to approximate the ILU-Surrogate accurately enough. 
We use \includegraphics[width=6pt]{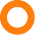} and \includegraphics[width=6pt]{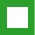} to depict these scenarios, where the degrees in x- and y-direction are fixed to 1 and 2, respectively, and only the degree in z-direction is variable.
In this case, the additional costs due to the higher degrees in z-direction are negligible and are not visible in both plots.
The runtime for our surrogate based forward-substitution is consistently lower than for the matrix-based version, even though we only get a minor performance gain.

Overall, we have presented a memory efficient ILU-Smoother which can be applied to much larger problems as if matrix-based implementations are used.
In the most general case, our method yields a modest runtime penalty for larger polynomial degrees. 
However, due to the reduced memory traffic, we can even achieve a small runtime gain with our surrogate approach for many relevant cases.

\section{Conclusion}\label{sec:5_conclusion}
In this paper, we have introduced an ILU algorithm on hybrid grid geometries and investigated its performance within a multigrid solver.
We replaced the matrix-based algorithm by two matrix-free variants based on surrogate polynomials approximating the stencils of the ILU matrix.
To our knowledge, this is the first matrix-free realization of a non-local operator based on an algebraic factorization.
Both matrix-free methods could attain the asymptotic convergence rates of their matrix-based counterparts.

The ILU convergence rates were robust for individual distorted tetrahedra which suggests a large performance gain for the matrix-free evaluation of the Steklov--Poincar\'e operator on distorted hybrid grids compared to simpler smoothing schemes as the GS algorithm.
The surrogate ILU provides a huge memory reduction with a small improvement in runtime compared to the usual matrix-based realizations.

Furthermore, only the standard ILU-Smoother was investigated.
The extension of our algorithm to a matrix-free version of the modified ILU described in Remark~\ref{remark:modifiedILU} and an investigation of its performance are still of interest.

\paragraph{Code availability}
The software used for obtaining the presented results are part of the open source framework HyTeG \cite{kohl2019hyteg} and publicly available at \cite{ilusource}. The results can be reproduced by executing the Python scripts in \url{apps/ILUSmoother/scripts}.

\appendix
\section{Smoother performance on single macro-tetrahedra}

\label{sec:appendix}
The following numerical experiments investigate the performance of our smoother on single macro-tetrahedra. 

\begin{figure*}[ht]
  \begin{center}
    \includegraphics[width=\textwidth]{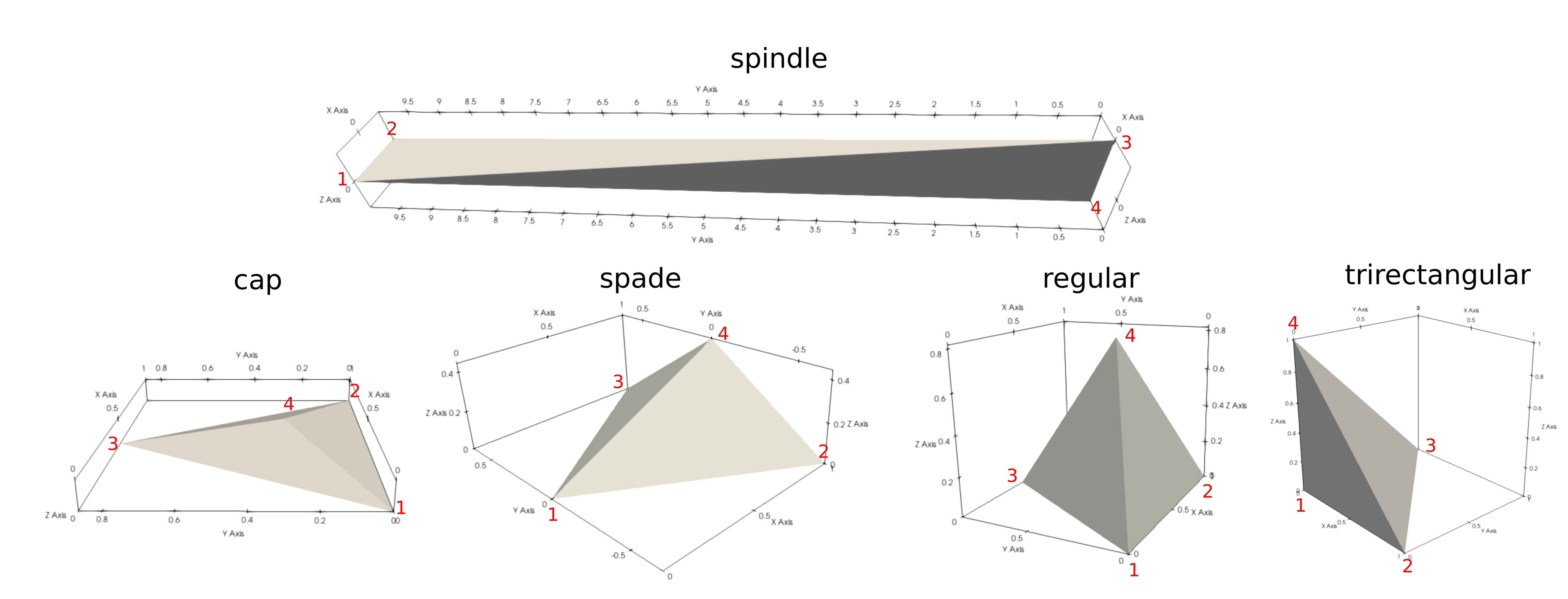}
  \end{center}
  \caption{\label{fig:regular-spade-cap-spindle}
    Tetrahedron shapes.
  }
\end{figure*}
In Figure~\ref{fig:regular-spade-cap-spindle}, we depict different tetrahedral shapes based on the ones in \cite{gmeiner2013optimization}, which are
the ``Spindle'',
the ``Cap'',
the ``Spade'',
and a regular tetrahedron.
The red numbers are used to assign numbers to the vertices to which a permutation will be applied.
The tetrahedral coordinates are given in Table~SM1 in the supplementary materials.

Due to the symmetry of the shapes, we only obtain a small number of different permutations.
The asymptotic convergence rates for these are given in Table~\ref{tab:asymptotic-convergence-rates-shapes} for both our ILU and symmetric GS-Smoother.
The rate determined by our permutation heuristic is marked in bold.
For all shapes, our ILU-Smoother yields better rates than the GS-Smoother if the best orientation is used.
This gain is especially dramatic for the ``Cap'' shape, where $\rho_{\textmd{GS}} \approx 0.5$, which means that 6-7 MG iterations are necessary with GS when compared to a single MG iteration with an ILU-Smoother.
In addition, the optimal permutations for the GS and the ILU-Smoother always coincide.
For a fair comparison, we will therefore use in all our future comparisons the same ordering for the GS as for the ILU-Smoother.
\begin{table}[ht]
  \begin{tabular}{l|ccc|cccc}
    $\pi$
     & $(1\,2\,3\,4)$               & $(1\,3\,2\,4)$ & $(1\,4\,2\,3)$                       & $(1\,2\,3\,4)$               & $(1\,2\,4\,3)$ & $(1\,3\,4\,2)$ & $(2\,3\,4\,1)$     \\
    \hline
    $\rho_{\textmd{GS}}$
     & 0.77                         & 0.54           & 0.78                                 & 0.52                         & 0.53           & 0.52           & 0.51               \\
    $\rho_{\textmd{ILU}}$
     & 0.65                         & {\bf 0.39}     & 0.35                                 & 0.010                        & 0.43           & 0.43           & {\bf 0.0096}        \\
    \hline
    $\#_{\textmd{GS}}$ & 53 & 23 & 56 & 22 & 22 & 22 & 21
    \\
    $\#_{\textmd{ILU}}$ & 33 & 15 & 14 & 3 & 17 & 17 & 3
    \\
    \hline
     & \multicolumn{3}{c|}{Spindle}
     & \multicolumn{4}{c}{Cap}
  \end{tabular}
  \vspace*{1cm} \\
  \begin{tabular}{l|cccccc|c}
    $\pi$
     & $(1\,2\,3\,4)$             & $(1\,2\,4\,3)$ & $(1\,3\,4\,2)$ & $(2\,1\,3\,4)$ & $(2\,1\,4\,3)$ & $(2\,3\,4\,1)$
     & $(1\,2\,3\,4)$
    \\
    \hline
    $\rho_{\textmd{GS}}$
     & 0.20                       & 0.085          & 0.20           & 0.079          & 0.20           & 0.055
     & 0.054
    \\
    $\rho_{\textmd{ILU}}$
     & 0.084                & 0.053          & 0.060          & {\bf 0.014}          & 0.14           & 0.028
     & {\bf 0.025}
    \\
    \hline
    $\#_{\textmd{GS}}$ & 9 & 6 & 9 & 6 & 9 & 5 & 5
    \\
    $\#_{\textmd{ILU}}$ & 6 & 5 & 5 & 4 & 8 & 4 & 4
    \\
    \hline
     & \multicolumn{6}{c|}{Spade}
     & Regular
  \end{tabular}
  \caption{\label{tab:asymptotic-convergence-rates-shapes}
    Asymptotic convergence rates $\rho$ and the number of iterations $\#$ to decrease the error by $10^{-6}$ for different permutations $\pi$.
    The convergence rates due to our reordering algorithm are marked in bold.
  }
\end{table}

\begin{figure*}[ht]
  \center
  \includegraphics[width=\textwidth]{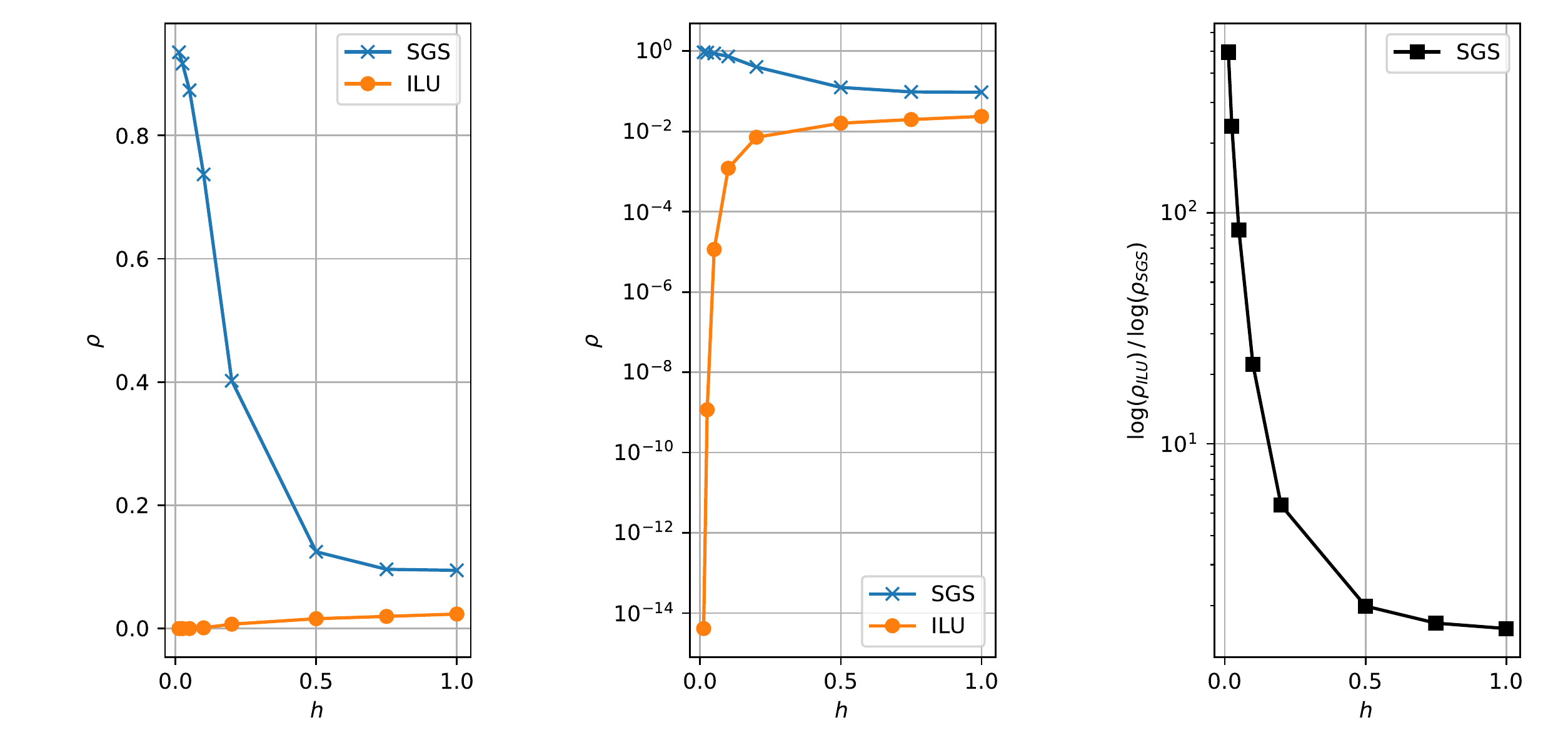}
  \caption{\label{fig:heights-single-tet}
    The multigrid algorithm on a distorted tetrahedron.
    Left: The asymptotic convergence rates for the multigrid algorithm.
    Middle: The convergence rates in a semilogarithmic plot.
    Right: Comparison of $\log(\rho_{ILU}) / \log(\rho_{SGS})$.
  }
\end{figure*}
We depict the convergence rate of a trirectangular tetrahedron for different heights of the top vertex in Figure~\ref{fig:heights-single-tet}.
As we see in the leftmost plot, decreasing the height degenerates the asymptotic convergence rate $\rho$ for multigrid with symmetric GS smoothing.
It is a well-known fact that the ILU-Smoother remains robust with respect to these deformations.
The logarithmic plot in the middle even suggests, that it becomes an exact solver.
In the right we have plotted $\log(\rho_{ILU}) / \log(\rho_{SGS})$, which is a measure how many additional multigrid iterations the SGS based multigrid need to achieve the same convergence rate as the ILU based multigrid. 
The performance gain due to an ILU becomes arbitrary large for small heights $h$.

Both Table~\ref{tab:asymptotic-convergence-rates-shapes} and Figure~\ref{fig:heights-single-tet} suggest, that the ILU on the reoriented tetrahedra mesh provides a robust smoother inside the macro-tetrahedrons. This transfers directly to the Steklov--Poincar\'e operator in which we use our ILU as an inner solver component.

\setcounter{figure}{0}
\setcounter{table}{0}
\setcounter{section}{0}
\setcounter{subsection}{0}
\renewcommand{\figurename}{Fig. SM}
\renewcommand{\tablename}{Table SM}
\renewcommand{\thesection}{SM\arabic{section}}
\renewcommand{\thesubsection}{SM\arabic{section}.\arabic{subsection}}
\section{Supplementary Material: Hybrid-Smoother}

In this second application, we extend the subgrid ILU-Smoother to a smoother on the global grid. 
On the interfaces between the structured grids, we apply a simple Gauss--Seidel-Smoother (GS-Smoother).
Thereby, we obtain a block smoother that we examine on a hybrid structured tetrahedral grid for a finite-element discretization of piecewise continuous finite-element functions.

Given a set of preconditioner matrices $C_{p} \in \R^{|\mcI_p| \times |\mcI_p|}$ for the lower dimensional primitives $p \in \mcV_H \cup \mcE_H \cup \mcF_H$
and a set of symmetric preconditioner matrices $C_{t,sym} = C^T_{t,sym} \in \R^{|\mcI_t| \times |\mcI_t|}$ for each tetrahedron $t \in \mcT_H$, we define a symmetric hybrid smoother in Algorithm~\ref{alg:hybrid-smoother}.
\begin{algorithm}[ht]
  \begin{algorithmic}[1]
  \State $ \bfu \gets \bfu + \sum_{v \in \mcV_H} R_v^T C^{-1}_v R_v\, (\bff - A \bfu) $. 
  \State $ \bfu \gets \bfu + \sum_{e \in \mcE_H} R_e^T C^{-1}_e R_e\; (\bff - A \bfu) $. 
	\State $ \bfu \gets \bfu + \sum_{f \in \mcF_H} R_f^T C^{-1}_f R_f (\bff - A \bfu) $. 
	\State $ \bfu \gets \bfu + \sum_{t \in \mcT_H} R_t^T C^{-1}_{t,sym} R_t (\bff - A \bfu) $. 
	\State $ \bfu \gets \bfu + \sum_{f \in \mcF_H} R_f^T C^{-T}_f R_f (\bff - A \bfu) $. 
	\State $ \bfu \gets \bfu + \sum_{e \in \mcE_H} R_e^T C^{-T}_e R_e\; (\bff - A \bfu) $. 
	\State $ \bfu \gets \bfu + \sum_{v \in \mcV_H} R_v^T C^{-T}_v R_v\, (\bff - A \bfu) $. 
  \end{algorithmic}
  \caption{Symmetric hybrid smoother.}
  \label{alg:hybrid-smoother}
\end{algorithm}
This smoother sequentially updates the degrees-of-freedom (DoF) on the macro-vertices, macro-edges, macro-faces and macro-tetrahedra.
In order to retain symmetry, we traverse the macro-hierarchy in the reverse direction by applying smoothing steps to the macro-face, macro-edge and finally to the macro-vertex DoF.
Note that procedures like that can be performed efficiently for distributed memory parallelizations since all DoF sharing a primitive type can be smoothed in parallel \cite{kohl2019hyteg}.
Due to the symmetrization it is possible to use our multigrid algorithm as a preconditioner for a conjugate gradient (CG)
or minimial residual (MINRES) method.

We define a hybrid symmetric Gauss--Seidel-Smoother (SGS-Smoother) as a reference example. 
It is given by
\[
	C_v = L_{A,v},\,
	C_e = L_{A,e},\,
	C_f = L_{A,f},
	\textmd{ and }
	C_{t,sym} = (L_{A,t}) (D_{A,t})^{-1} (L_{A,t})^T,
\]
for $v \in \mcV_H$, $e \in \mcE_H$, $f \in \mcF_H$ and $t \in \mcT_H$.
For simplicity, we will refer to it by Gauss--Seidel even though this is technically not correct since some of the DoF on the interfaces $\mcV_H \cup \mcE_H \cup \mcF_H$ are handled additively instead of multiplicatively.
Numerical experiments in \cite{kohl2019hyteg} have shown that in practice this has no impact on the convergence rates.

The hybrid ILU-Smoother is defined similarly by
\[
	C_v = L_{A,v},\,
	C_e = L_{A,e},\,
	C_f = L_{A,f},
	\quad \textmd{ and } \quad
	C_{t,sym} = L_t D_t L_t^T
  .
\]
In a nutshell, the smoother consists of Gauss--Seidel steps on the macro-vertices, macro-edges, and macro-faces. The ILU itself is just applied inside the macro-tetrahedra.

\begin{remark}
We stress that the submatrices 
$R_v A R^T_v$, $R_e A R^T_e$, and $R_f A R^T_f$ for $v \in \mcV_H, e \in \mcE_H,$ and $f \in \mcF_H$
on the lower dimensional primitives are well conditioned matrices.
Hence, even a simple iterative scheme is sufficient for reducing the error significantly.
Therefore, using an ILU also on the lower primitives is not necessary and only complicates the implementation.
\end{remark}

\subsection{Numerical experiments}
We investigate the performance of our hybrid ILU-Smoother on a hybrid mesh consisting of several tetrahedra.
For this, we will introduce a suitable benchmark problem. 
Firstly, we will directly apply the multigrid algorithm as solver and in a second scenario we will use the multigrid as a preconditioner within a CG method.

\paragraph{Benchmark}
Figure~\ref{fig:benchmark-scenario} depicts our benchmark scenario.
A unit cube with Dirichlet boundary conditions on the top and bottom and Neumann boundary conditions on the remaining sides
is divided into two parts where $h_{\textmd{lower}}$ is the height of the lower volume.
In both volumes, we assume a piecewise constant material parameter~$\kappa$ with a possible jump at the interface.
We want to compare the convergence rates for different heights of the lower cube.
This is similar to a scenario in geophysics in which the viscosity jumps from a narrow upper layer in the Earth's mantle that can barely be resolved to a wider lower layer \cite{davies1992mantle}.
We will use this scenario to assess the performance of our hybrid ILU-Smoother.

\begin{figure*}[ht]
  \begin{center}
    \includegraphics[width=.49\textwidth]{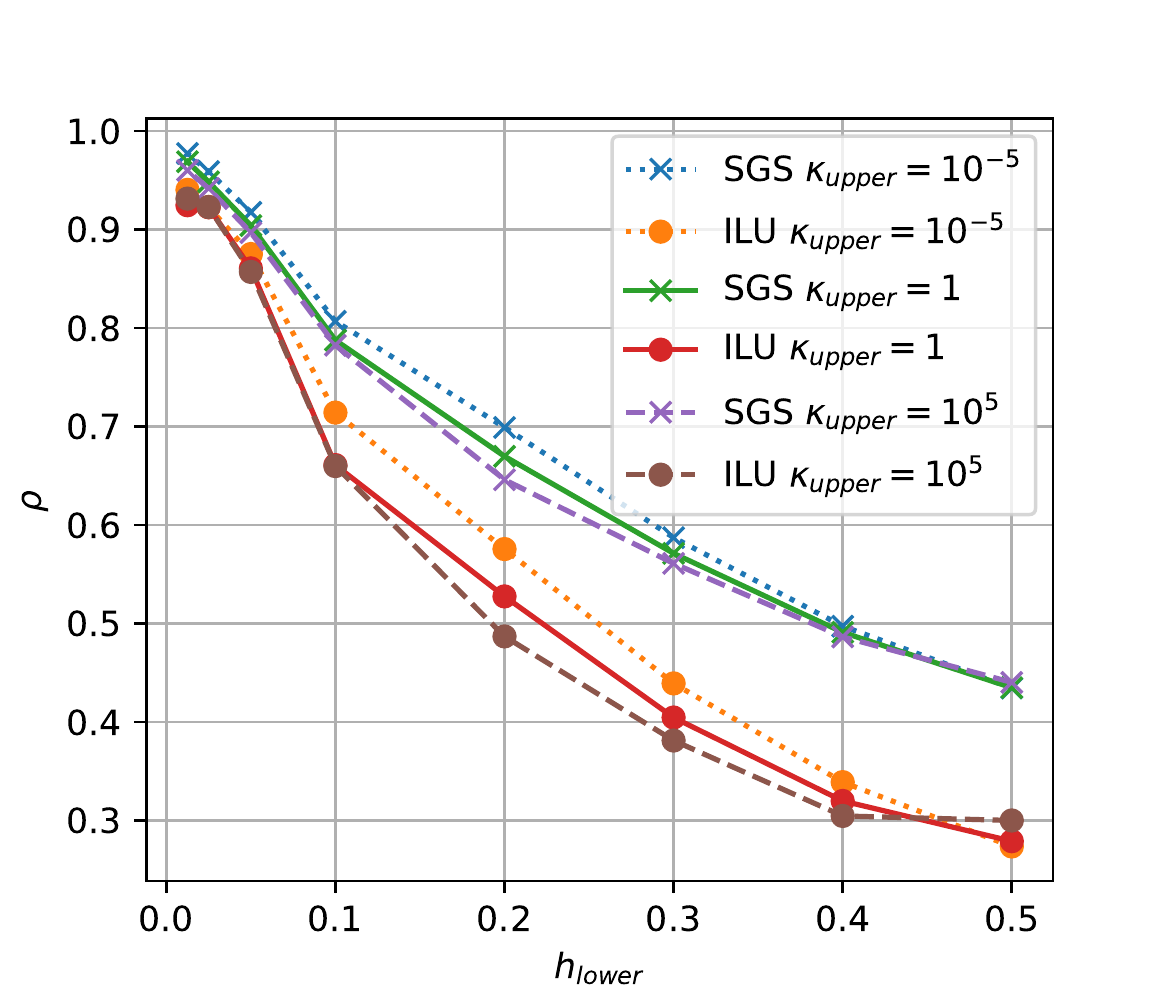}    \includegraphics[width=.49\textwidth]{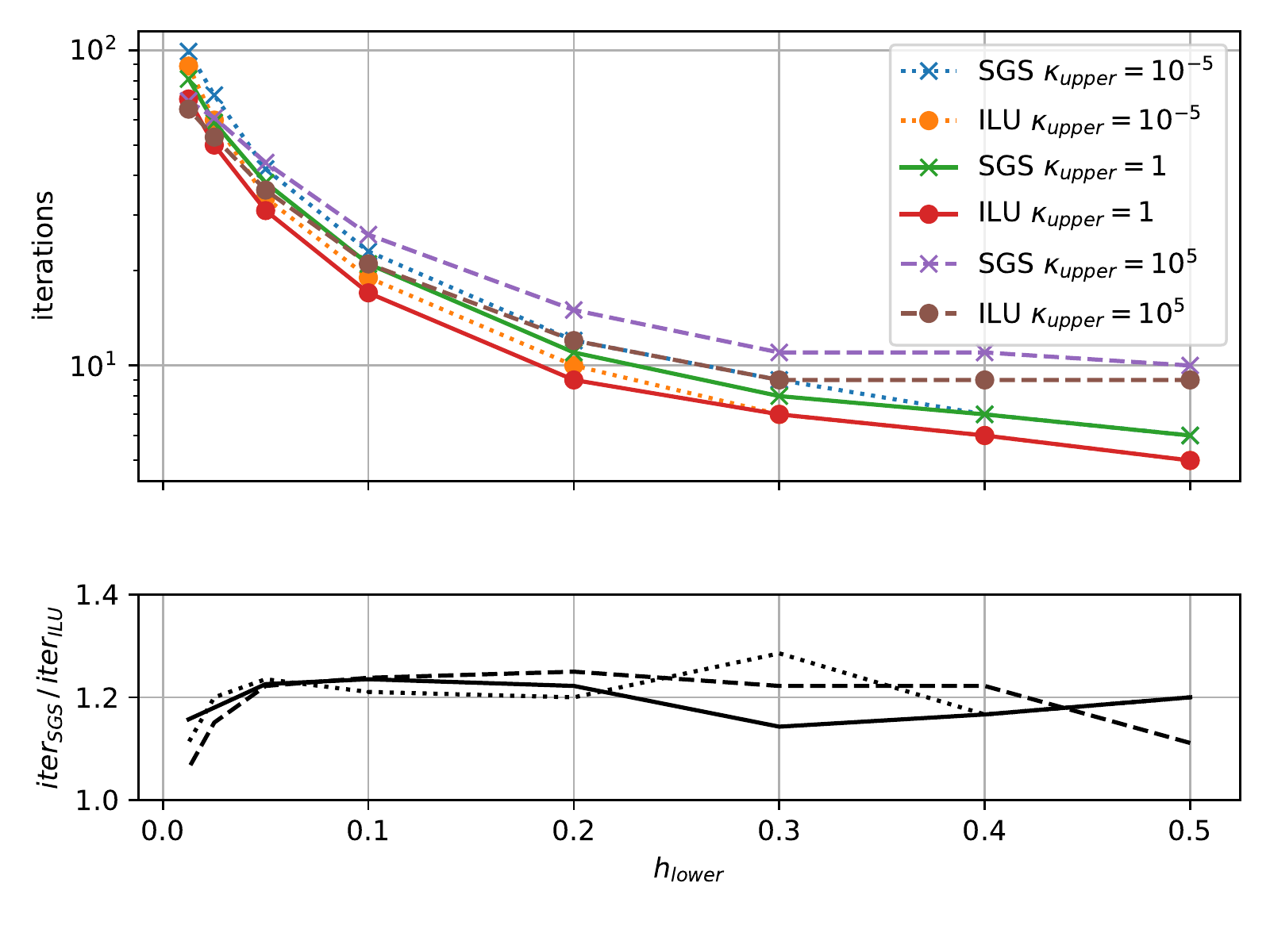}
  \end{center}
  \caption{\label{fig:hybridheight-mg-and-cg}
    Different solvers applied to our benchmark scenario (Fig.~\ref{fig:benchmark-scenario}).
    Left: Asymptotic convergence rates with a multigrid V-Cycle on grid levels 2 to 6, and 3 pre- and postsmoothing steps.
    Right: The V-Cycle used as a preconditioner for a PCG solver with iteration numbers in the upper plot and the ratio of the SGS and ILU iteration numbers in the lower plot.
    }
\end{figure*}

\begin{minipage}[ht]{0.4\textwidth}
  \begin{center}
    \vspace{1.25cm}
    \includegraphics[width=\textwidth]{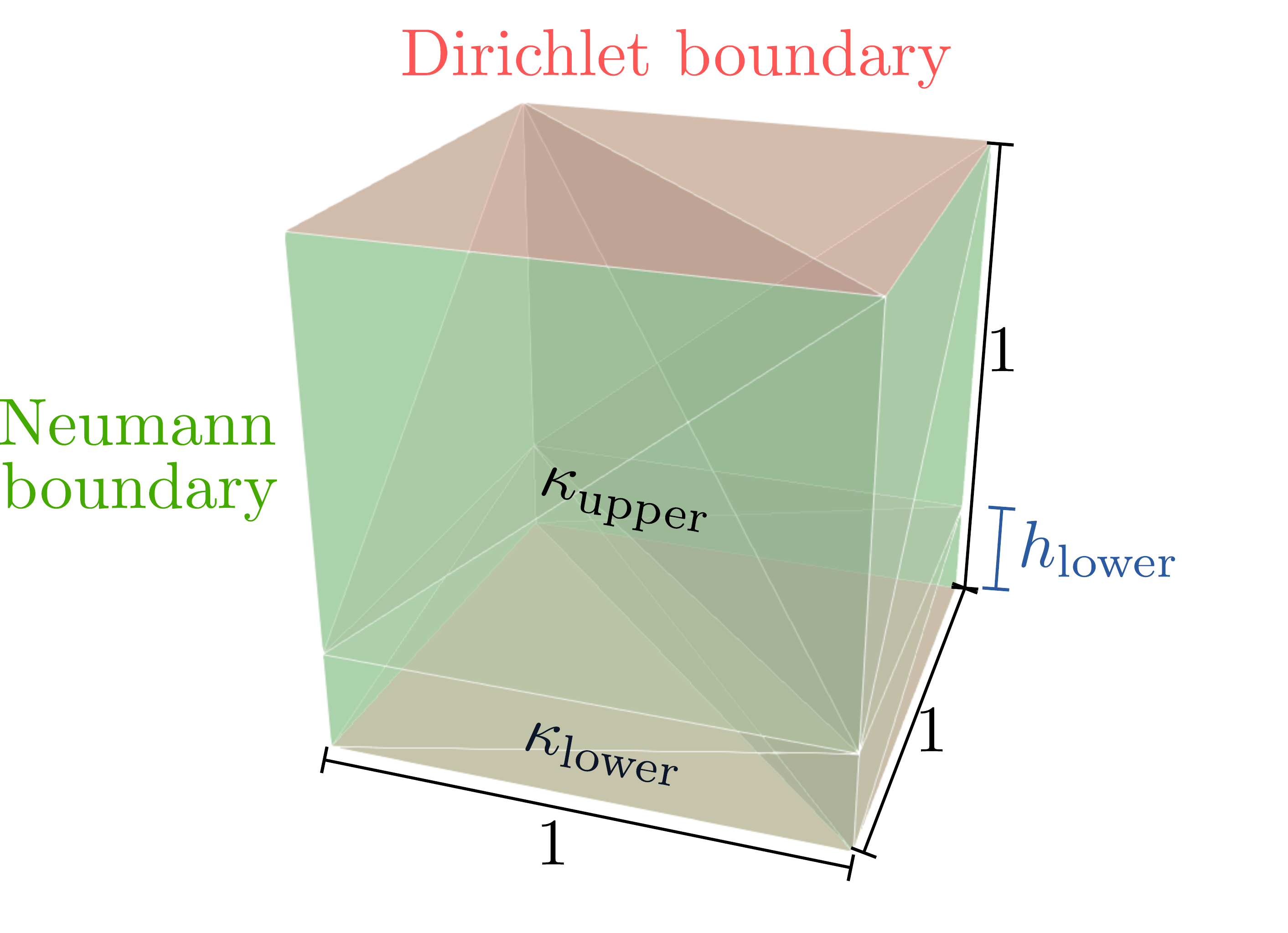}
    \vspace{0.90cm}\,
  \end{center}
  \captionof{figure}{\label{fig:benchmark-scenario}
    Benchmark scenario: The size of the lower half of the cube is decreased to approximate a tiny layer in the Earth's mantle.
    }
\end{minipage}\begin{minipage}[ht]{0.025\textwidth}\,
\end{minipage}\begin{minipage}[ht]{0.45\textwidth}
\centering
\begin{tabular}{|l|ccc|}
\hline
Geometry & x & y & z \\ \hline
Spindle                        & 0.0                    & 0.0                    & 0.5                    \\
                               & 0.0                    & 0.0                    & -0.5                   \\
                               & 0.5                    & 1.0                    & 0.0                    \\ 
                               & -0.5                   & 1.0                    & 0.0                    \\
                               \hline
Cap                            & 0.0                    & 0.0                    & 0.0                    \\
                               & 1.0                    & 0.0                    & 0.0                    \\
                               & 0.5                    & 0.866                  & 0.0                    \\
                               & 0.5                    & 0.288                  & 0.093                  \\
                               \hline
Spade                          & 0.0                    & 0.0                    & 0.0                    \\
                               & 1.0                    & -0.666                 & 0.0                    \\
                               & 1.0                    & 0.666                  & 0.0                    \\
                               & 1.0                    & 0.0                    & 0.443                 \\
                               \hline
\end{tabular}\captionof{table}{\label{tab:tetrahedral-coordinates} Tetrahedral coordinates for the Spindle, Cap and Spade.}
\end{minipage}
\paragraph{Multigrid solver}
The results for the asymptotic convergence rates of our multigrid algorithm are depicted in the left of Figure~\ref{fig:hybridheight-mg-and-cg},
for a fixed $\kappa_{lower} = 1$ and 3 different choices of $\kappa_{upper}$.
As expected, since the material jump happens on an interface that can be resolved by all coarse grids, it does not impact the SGS-Smoother (see \cite[Sec.~10.3]{hackbusch2013multi}) and the same is true for our ILU-Smoother.
For jumps inside the elements, special interpolation operators would have to be introduced \cite{alcouffe1981multi}.
Our ILU-Smoother consistently performs better than the SGS-Smoother for the multigrid algorithm.
For small tetrahedral heights, the performance of both smoothers decays and the convergence rate approaches one.
The slope of the performance degradation is similar for both smoothers.
Note that on the interfaces, we use the same Gauss--Seidel smoothing strategy for both smoothers which suggests that the performance loss is related to them and better methods are necessary to handle the interfaces.

\paragraph{Multigrid preconditioner}
In the right plot of Figure~\ref{fig:hybridheight-mg-and-cg}, we use our multigrid algorithm only as a preconditioner for a CG solver and plot the iteration count for getting the absolute unpreconditioned residual below $10^{-5}$.
The hybrid ILU-Smoother is an improvement compared to the SGS-Smoother, but again the hybrid solver cannot preserve the robustness with respect to degenerated tetrahedra.
The ratio between the number of iterations between a GS and an ILU based multigrid preconditioner is depicted in the lower right.
For a large range of heights it stays around 1.2 before it degenerates to 1 for extremely small heights.

\subsection{Conclusion}
We extended the smoother to a hybrid ILU-Smoother to directly use it as a smoother on a more complex domain with several tetrahedrons joined via an interface. 
In this case, our algorithm outperforms the SGS-Smoother while having comparable computational costs.

Although, a hybrid block smoother is perfectly suited for matrix-free approaches on hybrid meshes, we observed that it lacks robustness with respect to poorly structured macro-meshes.
Future work might therefore include applying more sophisticated techniques on the interfaces to further improve the performance.

\section{Simulation data}

Table~\ref{tab:tetrahedral-coordinates} provides the detailed coordinates of the reference tetrahedra for result reproductions and comparisons.

\bibliographystyle{siamplain}
\bibliography{references}
\end{document}